\documentclass[10pt,pdftex,a4paper]{amsart}
\usepackage{boxedminipage}
\usepackage{algpseudocode}
\usepackage{booktabs}
\usepackage{tabularx}
\usepackage{multirow}
\usepackage{subfigure}
\usepackage{mathrsfs} 
\renewenvironment{itemize}{\begin{list}{\labelitemi}{\leftmargin=1em} }
{\end{list}}
      
      \usepackage{graphicx}
      \usepackage{pdfsync}
      \usepackage{graphicx}
      \usepackage{enumerate}
      \usepackage{xspace}
      \usepackage{boxedminipage}
      \usepackage{algorithmicx}
      \usepackage[ruled]{algorithm}
      \usepackage{subfigure}
      \usepackage{algpseudocode}
      \usepackage[square,sort&compress,numbers]{natbib}
      \usepackage{hyperref} 
\usepackage{times}
      \usepackage{graphicx}
      \usepackage{epsfig}
      \usepackage{latexsym}
      \usepackage{amsmath,amssymb,amsxtra,amsthm}
      \usepackage{subfigure}
      \usepackage{stmaryrd}
      \usepackage{geometry}                
      \usepackage{commath}
      \usepackage{tikz}
      \usetikzlibrary{calc}
      
      \newcommand{\Proofname}{Proof}
      \newenvironment{Proof}[1][]
      	       {\par\noindent\textbf{\Proofname\xspace#1.\ }}
      	       {{\raggedright{{ }\hfill\qed}}} 
      \newtheoremstyle%
          {plain}
          {}
          {}
          {\mdseries\slshape}
          {}
          {\bfseries}
          {.}
          {1.0ex}
          {}
          
      \newtheoremstyle
          {note}
          {}
          {}
          {}
          {}
          {\bfseries}
          {.}
          {1.0ex}
          {}
          
      \swapnumbers
      \numberwithin{equation}{section}
      \newenvironment{The}[1][]{\subsection{Theorem\xspace{\ifx&#1&{}\else{ (#1)}\fi}}\slshape}{\upshape}
      \newenvironment{Prop}[1][]{\subsection{Proposition\xspace{\ifx&#1&{}\else{ (#1)}\fi}}\slshape}{\upshape}
      
      \newenvironment{Lem}[1][]{\subsection{Lemma\xspace{\ifx&#1&{}\else{ (#1)}\fi}}\slshape}{\upshape}
      \newenvironment{Cor}[1][]{\subsection{Corollary\xspace{\ifx&#1&{}\else{ (#1)}\fi}}\slshape}{\upshape}
      
      \theoremstyle{note}

      \newenvironment{Rem}[1][]{\subsection{Remark\xspace{\ifx&#1&{}\else{ (#1)}\fi}}}{}
      
      \newenvironment{Defn}[1][]{\subsection{Definition\xspace{\ifx&#1&{}\else{ (#1)}\fi}}}{}
      \newenvironment{Assum}[1][]{\subsection{Assumption\xspace{\ifx&#1&{}\else{ (#1)}\fi}}}{}
      \newenvironment{Pbm}[1][]{\subsection{Problem\xspace{\ifx&#1&{}\else{ (#1)}\fi}}}{}

      \newcommand{\Program}[1]{\textsf{#1}\xspace}
      \newcommand{\ALBERTA}{\Program{ALBERTA}}
      \newcommand{\PARAVIEW}{\Program{PARAVIEW}}
      \renewcommand{\vec}[1]{\ensuremath{\boldsymbol{#1}}}
      
      \newcommand{\pdt}{\ensuremath{\partial_t}}

      \newcommand{\Reals}{\ensuremath{{\mathbb{R}}}}
      \newcommand{\Integers}{\ensuremath{{\mathbb{Z}}}}
      
      \renewcommand{\O}{\ensuremath{{\Omega}}}
      \newcommand{\Ot}{\ensuremath{{\Omega_t}}}
      \newcommand{\Otn}{\ensuremath{{\Omega_{t^n}}}}
      
      \newcommand{\Oc}{\ensuremath{{\hat{\Omega}}}}
      \newcommand{\normal}{\ensuremath{{\vec{\nu}}}}
      
      \newcommand{\transpose}{\ensuremath{{\mathsf{T}}}}
      \newcommand{\Lp}[1]{\ensuremath{\operatorname{L}_{#1}}}
      \newcommand{\Sob}[2]{\ensuremath{\operatorname{W}^{#1,#2}}}
      \newcommand{\Hil}[1]{\ensuremath{\operatorname{H}^{#1}}}
      \newcommand{\V}{\ensuremath{{\mathbb{V}}}}
      \newcommand{\Vt}{\ensuremath{{\mathbb{V}_t}}}
      \newcommand{\Vtn}{\ensuremath{{\mathbb{V}^n}}}
      \newcommand{\Vc}{\ensuremath{\hat{\mathbb{V}}}}
      \newcommand{\ltwop}[3]{\ensuremath{\left\langle#1,#2\right\rangle}_{#3}}

      \newcommand{\seminorm}[3][]{\ensuremath{#1|#2#1|}_{#3}}
      \newcommand{\mynorm}[3][]{\ensuremath{#1\|#2#1\|}_{#3}}
      \newcommand{\ltwon}[3][]{\mynorm[#1]{#2}{\Lp{2}{#3}}}
      \newcommand{\linfn}[3][]{\mynorm[#1]{#2}{\Lp{\infty}{#3}}}
      \newcommand{\Hiln}[4][]{\mynorm[#1]{#2}{\Hil{#3}{#4}}}
      \newcommand{\A}{\ensuremath{\vec{\mathcal{A}}}}
      \newcommand{\lv}{\ensuremath{\left\vert}}
      \newcommand{\rv}{\ensuremath{\right\vert}}
      \usepackage{upgreek}
      \newcommand{\lap}{\ensuremath{\Updelta}}
      
      \newcommand{\T}{\ensuremath{{\mathscr{T}}}}
      \newcommand{\Tc}{\ensuremath{\hat{\mathscr{T}}}}
      \newcommand{\s}{\ensuremath{{{s}}}}

      \newcommand{\Th}{\ensuremath{\hat{\T}}}
      \newcommand{\Ritz}[1][t]{\ensuremath{{R}_{#1}}}
      \newcommand{\fim}{for $i=1,\dotsc,m$, }

      \newcommand{\fracl}[2]{#1/#2}
      \newcommand{\Clement}{\ensuremath{{\mathcal{I}^h}}}
      \newcommand{\Lagrange}{\ensuremath{{\Lambda^h}}}
      \newcommand{\Q}{\ensuremath{\mathcal{Q}}}
      \usepackage{color}

      \definecolor{MyGreen}{rgb}{0.05,0.4,0.05}
      \definecolor{RedViolet}{rgb}{0.1,0.1,0.75}

         \providecommand{\Cont}[1]{\ensuremath{\operatorname{C}^{#1}}}
         \usepackage{verbatim}
         \definecolor{SussexFlint}{rgb}{.00,.19,.21}
         \definecolor{SussexGrey}{rgb}{.51,.58,.49}
         \definecolor{SussexOrange}{rgb}{.94,.29,.00}
         \definecolor{SussexYellow}{rgb}{1.00,.73,.00}
         \definecolor{SussexRed}{rgb}{.94,.01,.49}
         \definecolor{SussexPurple}{rgb}{.48,.06,.44}
         \definecolor{SussexGreen}{rgb}{.00,.58,.46}
         \definecolor{SussexBlue}{rgb}{.00,.58,.65}
         \colorlet{a}{SussexOrange}
         \colorlet{b}{SussexRed}
         \colorlet{c}{SussexYellow}
         \colorlet{d}{SussexPurple}
         \colorlet{e}{SussexGreen}
         \colorlet{f}{SussexBlue}
         \colorlet{g}{white}
         \colorlet{h}{SussexGrey}
         \colorlet{i}{black}
         \colorlet{j}{SussexFlint}
         \newcommand{\bbil}[4][t]{\ensuremath{b_{#1}\left({#2},{#3}\right)}}
         
        \usepackage{ifthen}
        \usepackage{amsmath}
        \usepackage{amssymb}    
        \usepackage{stmaryrd} 
        \usepackage{eufrak}     
        \usepackage{mathrsfs} 
        \usepackage{bbold}    
        \provideboolean{usemathrsfs}
        \setboolean{usemathrsfs}{true}
        \usepackage{enumerate}
        \usepackage{xspace}
        \usepackage[greek,british,english]{babel}
        \usepackage{verbatim}
        \usepackage{fancyvrb}
        \usepackage{hyperref}
        \usepackage{listings}
        \usepackage{xcolor}
        \provideboolean{usebeamer}
        \provideboolean{isthesis}
        \provideboolean{isamsltex}
        \provideboolean{issiamltex}
    \usepackage{xcolor}
    \colorlet{a}{red}
    \colorlet{b}{blue}
    \colorlet{c}{green!50!blue}
    \colorlet{d}{magenta}
    \colorlet{e}{cyan}
    \colorlet{f}{yellow!50!black}
    \colorlet{g}{white}
    \colorlet{h}{black!50}
    \colorlet{i}{black}
    \colorlet{j}{black!75} 
    \providecommand{\linkedurl}[1]{\url{#1}}
    \providecommand{\linkedemail}[1]{\href{mailto:#1}{#1}}
    
    \providecommand{\email}[1]{{\linkedemail{#1}}}
    \providecommand{\Ignore}[1]{}
    \providecommand{\ignore}[1]{}
    \providecommand{\freeze}[1]{}
    \providecommand{\crossout}[1]{{\textcolor{red!20}{#1}}}
    \providecommand{\highlight}[1]{{\color{blue}#1}}
    \providecommand{\standout}[1]{\colorbox{a}{\textcolor{g}{#1}}}
    
    \provideboolean{shownotes}
    \setboolean{shownotes}{true}
    \newcounter{margnote}[page]
    \providecommand{\margnotemark}{{\standout{\upshape\texttt{\arabic{margnote}}}}}
    \providecommand{\margnote}[2][]{
      \ifthenelse{
        \boolean{shownotes}
      }{
        \stepcounter{margnote}
        \margnotemark
        \marginpar{
          \texttt{
            \begin{minipage}{2cm}
              \raggedright\tiny
              \margnotemark
              {\ifx&#1&{}\else{#1 says:}\fi} 
              #2
            \end{minipage}
          }
        }
      }{
      }
    }
    \providecommand{\mathnote}[2][]{
      \ifthenelse{
        \boolean{shownotes}
      }{
        \stepcounter{margnote}
        \text{
          \standout{
            \texttt{
              \tiny
              \margnotemark
              #1:
              #2
            }
          }
        }
      }{
      }
    }
    \provideboolean{showtodo}

    \providecommand{\todo}[1]{\ifthenelse{\boolean{showtodo}}{\margnote[To do.]{#1}}{}}
    \providecommand{\Todo}[1]{
      \ifthenelse{\boolean{showtodo}}{
        \begin{center}
        \begin{tikzpicture}
         \node[fill=a!17]{
           \begin{minipage}{\textwidth}
             \texttt{\bfseries{\small #1}}
           \end{minipage}
         };
        \end{tikzpicture}
        \end{center}
      }{}}

    \provideboolean{showcomments}
    \providecommand{\margincomment}[1]{
    \ifthenelse{\boolean{showcomments}}{\marginpar{\tiny #1}}{}
    }
    \provideboolean{showchanges}
    \setboolean{showchanges}{false}
    \providecommand{\changes}[1]{
      \ifthenelse{\boolean{showchanges}}{{\highlight{#1}}}{#1}
    }
    \providecommand{\changefromto}[3][replace with]{
      \ifthenelse{\boolean{showchanges}}
      {{\crossout{#2}\margnote{#1}}{\highlight{#3}}}
      {#3\xspace}
    }
    \providecommand{\ChangePar}[2]{
      \ifthenelse{\boolean{showchanges}}
      {{\par$\mapsfrom$ \textcolor{red!20}{#1}}{\par$\mapsto$ \textcolor{blue}{#2}}}
      {\par #2}
    }
    \providecommand{\InsertPar}[1]{
      \ifthenelse{\boolean{showchanges}}
      {{\par$\mapsto$ \textcolor{blue}{#1}}}
      {\par #1}
    }

    \usepackage{mathrsfs}
    \provideboolean{usemathrsfs}
    \setboolean{usemathrsfs}{true}
    \providecommand{\mathscript}
    	   {\mathscr}

     \providecommand{\bbbold}{\mathbb}

     \providecommand{\rR}{\ensuremath{\bbbold R}\xspace}

     \providecommand{\reals}{\rR}
     \providecommand{\R}[1]{\reals^{#1}}

     \providecommand{\Th}{\ensuremath{\varTheta}\xspace}

     \providecommand{\W}{\ensuremath{\varOmega}\xspace}

     \providecommand{\qp}[1]{\ensuremath{\left({#1}\right)}}
     \providecommand{\qpreg}[1]{\ensuremath{(#1)}}
     \providecommand{\qpbig}[1]{\ensuremath{\big(\!#1\!\big)}}
     \providecommand{\qpBig}[1]{\ensuremath{\Big(\!#1\!\Big)}}
     \providecommand{\qpbigg}[1]{\ensuremath{\bigg(\!\!#1\!\!\bigg)}}
     \providecommand{\qpBigg}[1]{\ensuremath{\Bigg(\!\!#1\!\!\Bigg)}}
     \providecommand{\qb}[1]{\ensuremath{\!\left[{#1}\right]}}

     \providecommand{\qc}[1]{\ensuremath{\left\{{#1}\right\}}}

     \providecommand{\clinter}[2]{\ensuremath{\left[#1,#2\right]}\xspace}

     \providecommand{\powqp}[2]{\ensuremath{\qp{#2}^{\kern -.2em\lower .2ex\hbox{\scriptsize $#1$}}\kern-.1em}}
     \providecommand{\powqpreg}[2]{\ensuremath{\qpreg{#2}^{\kern -.2em\lower .3ex\hbox{\scriptsize $#1$}}\kern-.3em}}
     \providecommand{\powqpbig}[2]{\ensuremath{\qpbig{#2}^{\kern -.2em\lower .3ex\hbox{\scriptsize $#1$}}\kern-.3em}}
     \providecommand{\powqpBig}[2]{\ensuremath{\qpBig{#2}^{\kern -.2em\lower .3ex\hbox{\scriptsize $#1$}}\kern-.3em}}
     \providecommand{\powqpbigg}[2]{\ensuremath{\qpbigg{#2}^{\kern -.2em\lower .3ex\hbox{\scriptsize $#1$}}\kern-.3em}}
     \providecommand{\powqpBigg}[2]{\ensuremath{\qpBigg{#2}^{\kern -.2em\lower .3ex\hbox{\scriptsize $#1$}}\kern-.3em}}

     \providecommand{\norm}[2][]{\ifx&#1&\left|\else#1|\fi#2\ifx&#1&\right|\else#1|\fi}
     
     \providecommand{\abs}[1]{\ensuremath{\left|#1\right|}}
     \providecommand{\Norm}[2][]{\ifx&#1&\left\|\else#1\|\fi#2\ifx&#1&\right\|\else#1\|\fi}

     \providecommand{\ltwop}[2]{\ensuremath{\qa{#1,#2}}}

     \providecommand{\setof}[1]{{\qc{#1}}}
     
     \providecommand{\seqof}[1]{\qp{#1}}
     \providecommand{\seq}[1]{\seqof{#1}}
     
     \providecommand{\seqi}[3][]{\seq{#2_{#3}}_{\ifx&#1&{#3}\else{#1}\fi}}
     \providecommand{\seti}[3][]{\setof{#2_{#3}}_{\ifx&#1&{#3}\else{#1}\fi}}
     \providecommand{\sequ}[3][]{\seq{#2^{#3}}_{\ifx&#1&{#3}\else{#1}\fi}}
     \providecommand{\setu}[3][]{\setof{#2^{#3}}_{\ifx&#1&{#3}\else{#1}\fi}}
     
     \providecommand{\setsi}[3]{\seti[{#2\in#3}]{#1}{#2}}

     \providecommand{\seqsufromto}[4]{\sequ[{\rangefromto{#2}{#3}{#4}}]{#1}{#2}}

     \providecommand{\integerbetween}[2]{\ensuremath{={#1},\dotsc,{#2}}}
     \providecommand{\rangefromto}[3]{\ensuremath{#1\integerbetween{#2}{#3}}}
    
    \providecommand{\registered}%
    {\ensuremath{^\text{\textregistered}}}
    
    \providecommand{\constant}[1]{\ensuremath{C_{#1}}}
    \providecommand{\constref}[2][]{\ensuremath{\constant{\text{\ref{#2}{\ifx&#1&{}\else{,\ensuremath{#1}}\fi}}}}}

    \providecommand{\fracl}[3][]{{#2}#1/{#3}}

    \providecommand{\Eye}[1]{
      \begin{bmatrix}
      \ifthenelse{#1>1}{
        \ifthenelse{#1>2}{
          \ifthenelse{#1>3}{
            1&0&\dotso&0
            \\
            0&1&\dotso&0
            \\
            \vdots&\vdots&\ddots&\vdots
            \\
            0&0&\dotso&1
          }{
            1&0&0
            \\
            0&1&0
            \\
            0&0&1
          }
        }{
          1&0
          \\
          0&1
        }
      }{
        1
      }
      \end{bmatrix}
    }
    
    \providecommand{\pd}[1]{\ensuremath{\partial_{#1}}\xspace} 
    
    \providecommand{\pdt}[1][]{\pd t{{\ifx&#1&{}\else{\qb{#1}}\fi}}}                       
    
    \providecommand{\transpose}{{\boldsymbol\intercal}}   

    \providecommand{\lap}{\ensuremath{\Delta}}
    \providecommand{\normal}{\vec n}

    \renewcommand{\vec}[1]{\ensuremath{\boldsymbol{#1}}}

    \providecommand{\Symmatrices}[1]{\ensuremath{\operatorname{Sym}{(\R{d\times d})}}}

    \providecommand{\Forall}{\:\forall\:}
    
    \providecommand{\Foreach}{\quad\Forall}

    \usepackage{amscd}

    \providecommand{\Program}[1]{\textsf{#1}\xspace}

    \providecommand{\euro}{\textgreek{\euro}}
    
    \providecommand{\ListParameters}{}
    \renewcommand{\ListParameters}
    {
    	 \setlength{\topsep}{0em}
    	 \setlength{\leftmargin}{0em}
             \setlength{\itemsep}{0ex}
    	 \setlength{\parsep}{.5ex}
    	 \setlength{\itemindent}{\labelsep}
    	 \addtolength{\itemindent}{\labelwidth}
    }
    
    \newcounter{LetterListItem}
    \renewcommand{\theLetterListItem}{(\alph{LetterListItem})}
    {
    	\begin{list}%
    	{\theLetterListItem\ }%
    	{\usecounter{LetterListItem}
    	 \ListParameters
    	}
    }%
    {\end{list}}
    \newcounter{NumberListItem}
    \renewcommand{\theNumberListItem}{\arabic{NumberListItem}}
    {
    	\begin{list}%
    	{\theNumberListItem.\ }%
    	{\usecounter{NumberListItem}%
    	 \ListParameters
    	}
    }%
    {\end{list}}
    \newcounter{QuestionListItem}
    \renewcommand{\theQuestionListItem}{\textbf{Question \arabic{QuestionListItem}}}
    {
    	\begin{list}%
    	{\theQuestionListItem.\ }%
    	{\usecounter{QuestionListItem}%
    	 \ListParameters
    	}
    }%
    {\end{list}}
    \newcounter{RomanListItem}
    \renewcommand{\theRomanListItem}{(\roman{RomanListItem})}
    {
    	\begin{list}%
    	{\theRomanListItem\ }%
    	{\usecounter{RomanListItem}
    	 \ListParameters
    	}
    }%
    {\end{list}}
    \newcounter{StepsItem}
    {
    	\begin{list}%
    	{Step \theStepsItem.\ }%
    	{\usecounter{StepsItem}%
    	 \ListParameters
    	}
    }%
    {\end{list}}
    
    \usepackage{amsthm} 
    \ifthenelse{\boolean{isthesis}}{%
      \setcounter{secnumdepth}{1}%
    }{
      \setcounter{secnumdepth}{2} 
    }
    \providecommand{\ListParameters}{}
    \renewcommand{\ListParameters}
    {
    	 \setlength{\topsep}{0em}
    	 \setlength{\leftmargin}{0em}
             \setlength{\itemsep}{0ex}
    	 \setlength{\parsep}{.5ex}
    	 \setlength{\itemindent}{\labelsep}
    	 \addtolength{\itemindent}{\labelwidth}
    }
    \newtheoremstyle{plain}
      {}
      {}
      {\mdseries\slshape}
      {\parindent}
      {\bfseries}
      {.}
      {.5em}
      {}
    
    \newtheoremstyle{note}
      {}
      {}
      {}
      {\parindent}
      {\bfseries}
      {.}
      {.5em}
      {}
    
    \newtheoremstyle{claim}
      {}
      {}
      {\mdseries\slshape}
      {}
      {\bfseries}
      {}
      {.5em}
      {}
    
    \newtheoremstyle{exercise}
      {}
      {}
      {}
      {}
      {\bfseries}
      {.}
      {1em}
      {}
    
    \newtheoremstyle{break}
      {}
      {}
      {}
      {}
      {\bfseries}
      {.}
      {\newline}
      {}
    
      \providecommand{\SolName}{Solution}
      \providecommand{\Proofname}{Proof}
    
    \ifthenelse{\boolean{isthesis}}{%
      
    }{%
      
    }
    \providecommand{\pdfformat}[1]{
       \provideboolean{pdfoutput}
       \setboolean{pdfoutput}{#1}
      \ifthenelse{\boolean{pdfoutput}}{
        \typeout{using pdf}
        \providecommand{\graphext}{pdf}
        \renewcommand{\graphext}{pdf}
        \providecommand{\graphextex}{pdf_t}
        \renewcommand{\graphextex}{pdf_t}
      }{
        \typeout{using eps}
        \usepackage[dvips]{graphicx,xcolor}
        \providecommand{\graphext}{eps}
        \renewcommand{\graphext}{eps}
        \providecommand{\graphextex}{eps_t}
        \renewcommand{\graphextex}{eps_t}
      }
      \usepackage{epsfig}
      \usepackage{tikz}
      \usepackage{rotating}
      \definecolor{SussexFlint}{rgb}{.00,.19,.21}
      \definecolor{SussexGrey}{rgb}{.51,.58,.49}
      \definecolor{SussexOrange}{rgb}{.94,.29,.00}
      \definecolor{SussexYellow}{rgb}{1.00,.73,.00}
      \definecolor{SussexRed}{rgb}{.94,.01,.49}
      \definecolor{SussexPurple}{rgb}{.48,.06,.44}
      \definecolor{SussexGreen}{rgb}{.00,.58,.46}
      \definecolor{SussexBlue}{rgb}{.00,.58,.65}
      \colorlet{a}{SussexOrange}
      \colorlet{b}{SussexRed}
      \colorlet{c}{SussexYellow}
      \colorlet{d}{SussexPurple}
      \colorlet{e}{SussexGreen}
      \colorlet{f}{SussexBlue}
      \colorlet{g}{white}
      \colorlet{h}{SussexGrey}
      \colorlet{i}{black}
      \colorlet{j}{SussexFlint}
      \newcommand{\mausDarkColorTheme}{
        \colorlet{a}{SussexYellow!50!yellow}
        \colorlet{b}{SussexGreen!50!green}
        \colorlet{c}{SussexBlue}
        \colorlet{d}{SussexOrange!50!yellow}
        \colorlet{e}{SussexRed!50!red}
        \colorlet{f}{SussexPurple!50!magenta}
        \colorlet{g}{black}
        \colorlet{h}{SussexFlint!50!black}
        \colorlet{i}{white}
        \colorlet{j}{SussexGrey}
      }
    }
    \providecommand{\solution}{\textbf{\SolName.}\xspace}
    
     \newcounter{phantombox}[enumi]
     \provideboolean{showphantoms}
     \renewcommand{\thephantombox}{\roman{phantombox}}
     \newcommand{\phantombox}[1]{\stepcounter{phantombox}
       \ensuremath{\boxed{
           \ifthenelse
    	   {\boolean{showphantoms}}
    	   {#1^{\phantom{\textup{(\thephantombox)}}}}
    	   {\phantom{#1}^{\textup{(\thephantombox)}}}
         }
       }
     }

     \provideboolean{hidesolution}
     \newcommand{\consolution}[1]{
       \ifthenelse{
         \boolean{hidesolution}
       }{
       }{
         {\par \small {\solution}\ #1\par\ \\[5pt]}}
     }
     \provideboolean{showmarks}
     \renewcommand{\marks}[1]{
       \ifthenelse{\boolean{showmarks}}{\marginpar{{\tiny [$#1$ marks]}}}{}}
     \newcommand{\condibreak}{\ifthenelse{\boolean{hidesolution}}{\newpage}{}}
     
\setboolean{shownotes}{true}
\setboolean{showchanges}{false}
\provideboolean{showcondense}
\setboolean{showcondense}{false}
\providecommand{\condense}[1]{
  \ifthenelse{\boolean{showcondense}}{{\color{red}{#1}}}{}
}
\provideboolean{showchangesvtwo}
\setboolean{showchangesvtwo}{false}
\providecommand{\changesvtwo}[1]{
  \ifthenelse{\boolean{showchangesvtwo}}{{\color{blue}{#1}}}{#1}
}

\provideboolean{showreviewadditions}
\setboolean{showreviewadditions}{false}
\providecommand{\additions}[1]{
  \ifthenelse{\boolean{showreviewadditions}}{{\color{blue}{#1}}}{#1}
}

\provideboolean{showreviewdeletions}
\setboolean{showreviewdeletions}{false}
\providecommand{\deletions}[1]{
  \ifthenelse{\boolean{showreviewdeletions}}{{\color{red}{#1}}}{}
}

\title[Imex and finite elements for reaction--diffusion systems on evolving domains]{
  Implicit--explicit timestepping with finite element approximation of reaction--diffusion systems on evolving domains
}
\date{\today}
\author[O.~Lakkis]{Omar Lakkis}
\address{
  Omar Lakkis,
  Department of Mathematics, University of Sussex,
  Brighton, England, GB-BN1~9QH.
}
\email{\linkedemail{o.lakkis@sussex.ac.uk}}
\urladdr{\url{http://www.maths.sussex.ac.uk/~omar}}
\author[A.~Madzvamuse]{Anotida Madzvamuse}
\address{
  Anotida Madzvamuse,
  Department of Mathematics, University of Sussex,
  Brighton, England, GB-BN1~9QH. 
}
\email{\linkedemail{a.madzvamuse@sussex.ac.uk}}
\author[C. Venkataraman]{Chandrasekhar Venkataraman}
\address{
  Chandrasekhar Venkataraman,
  {Department of Mathematics, University of Sussex,
  Brighton, England, GB-BN1~9QH. }
}
\email{{\linkedemail{c.venkataraman@sussex.ac.uk}}}
\thanks{
  AM would like to acknowledge partial financial
  support from the British Council through its \emph{UK-US New
    Partnership Fund (PMI2)}, the London Mathematical Society 
  ({\emph{R4P2}}) and the EPSRC small grant scheme: EP/H020349/1.
  \\
  The research of CV has been supported by the British Engineering
  and Physical Sciences Research Council (EPSRC), Grant EP/G010404.
  \\
  All authors acknowledge the EPSRC DTA Fellowship.
}
\begin{document}
\maketitle
\begin{abstract}
\additions{We present and analyse an implicit--explicit timestepping procedure
with finite element spatial approximation for a semilinear
reaction--diffusion systems on evolving domains arising from
biological models, such as Schnakenberg's (1979).  We employ a
Lagrangian formulation of the model equations which permits the error
analysis for parabolic equations on a fixed domain but introduces
technical difficulties, foremost the space-time dependent conductivity
and diffusion.  We prove optimal-order error estimates in the
$\Lp{\infty}(0,T;\Lp{2}(\W))$ and $\Lp{2}(0,T;\Hil{1}(\W))$ norms, and
a pointwise stability result.  We remark that these apply to Eulerian
solutions. Details on the implementation of the Lagrangian and the
Eulerian scheme are provided. We also report on a numerical
experiment for an application to pattern formation on an evolving
domain.} 
\end{abstract}
\section{Introduction}
\label{sec:intro}
  Since the seminal paper of \citet{turing1952cbm},
  \emph{time-dependent reaction--diffusion systems (RDSs)} have been
  studied as models for pattern formation in natural-process driven
  morphogenesis and developmental biology (see
  \citet{murray2003mathematical} for details).  An important
  generalisation of these models consists in considering RDSs posed on
  \emph{evolving domains}.  This stems from the now relatively
  well-known observation that in many cases growth of organisms plays
  a pivotal role in the emergence of patterns and their evolution
  during growth development \citep{murray2003mathematical,
    kondo1995rdw}.  RDSs on evolving domains have a wider scope of
  application, e.g., competing-species of micro-organisms in
  environmental biology, chemistry of materials and corrosion
  processes, the spread of pollutants. Numerical simulations of RDSs
  on time-evolving domains reproducing the empirically observed
  pattern formation processes are commonly used \citep{kondo1995rdw,
    barrio2009modeling, miura2006mixed, raquelsfem,amago,
    venk11chemotaxis}. It is essential for scientists to
  computationally approximate and appreciate the error between
  simulations and exact solutions of such RDSs. Galerkin finite
  elements~\citep{Thomee} are among the methods of choice to
  approximate such systems.

  \par In spite of their widespread use, to the best of our knowledge,
  no complete error analysis of approximating finite element schemes
  for \emph{nonlinear reaction--diffusion systems on evolving domains}
  is available in the literature, thus motivating this work.  This is
  a sibling paper to \cite{venkataraman2010global} where we analysed
  the well-posed nature of (exact) RDSs on evolving domains.  In most
  practical applications, the evolving domain is usually a surface
  embedded in the three-dimensional Euclidean space, but for
  simplicity we restrict our discussion to the case where both the
  reference domain and the evolving domain are flat, deferring thus
  the analysis of RDSs on evolving curved surfaces.

\par
\begin{Pbm}[RDS on a time-dependent evolving domain]
  \label{Pbm:RDS_td}
  We study a RDS, also considered in \citep{crampin1999, ano2000},
  which models a system of chemicals that interact through the reaction
  terms only and diffuse in the domain independently of each other. Given
  an integer $m\geq1$, the vector
  $\vec{u}\left(\vec{x},t\right)\in\mathbb R^m$, denoting the
  concentration of the chemical species $i=1,\dotsc,m$, at a spatial
  point $\vec{x} \in\Ot\subset \mathbb{R}^d$, $d=1,2,3$, at time $t
  \in [0,T]$, $T>0$,
  satisfies the following \emph{initial--boundary value problem}
  \begin{equation}\label{eqn:model_problem}
    \begin{gathered}
      \pdt{u}_i(\vec{x},t)-D_i\lap{u}_i(\vec{x},t)+\nabla\cdot\left[\vec{a}u_i\right](\vec{x},t)=f_{i}\left(\vec{u}(\vec{x},t)\right), 
      \quad
      \vec{x}\in\Ot,t\in(0,T],
        \\
          [\normal\cdot\nabla{u}_i](\vec{x},t)=0, 
          \quad
          \vec{x}\in\partial\Ot, t>0,
          \\
          u_i(\vec{x},0)=u_i^{0}(\vec{x}), 
          \quad
          \vec{x}\in\O_0,
   \end{gathered}
    \end{equation}
    where $\Omega_{t}$, detailed in \S\ref{sec:evolving-domain}, is a
    simply connected Lipschitz continuously \emph{evolving domain}
    with respect to $t\in[0,T]$, and
    $\vec{D}:=({D}_1,\dotsc,{D}_m)^\transpose$ is a vector of strictly
    positive diffusion coefficients.  Detailed assumptions on the
    \emph{nonlinear reaction vector field}
    $\vec{f}:=({f}_1,\dotsc,{f}_m)^\transpose$ are given in
    \S\ref{assum:reaction}.  The convection $\vec{a}=(a_1,\dotsc,a_d)^\transpose$ is induced by the material
    deformation due to the evolution of the domain.  The initial data
    $\vec{u}^{0}$ is a positive-entry bounded field.  Since we are
    primarily interested in pattern formation phenomena that arise as
    a result of self-organisation within a domain without
    outside-world communication we consider homogeneous Neumann
    boundary conditions, but other types of boundary conditions could
    be studied as well within our framework.
\end{Pbm}
\subsection{Main results}
  The core result in this paper is Theorem~\ref{thm:fully_discrete}
  where we prove optimal convergence rates of the discrete solution in
  $\Lp\infty(0,T;\Lp{2}(\Oc))^m$ and $\Lp2(0,T;\Hil1(\Oc)^m)$ (where
  $\Oc$ is a transformed version of $\Ot$ to be described next).  Our
  theoretical results are illustrated by numerical experiments,
  aimed mainly at quantifying the pattern formation phenomena related
  to the type of growth in the domains.
\subsection{A Lagrangian approach}
  We employ, both for the analysis and the implementation of the
  computational method, a Lagrangian formulation of Problem
  \ref{Pbm:RDS_td}, in the sense employed in fluid-dynamics, i.e.,
  where the evolving domain, $\Ot\in\mathbb R^d$, is the image of a
  time-dependent family of diffeomorphisms $\vec\A_t$ on a
  \emph{reference domain} $\Oc\in\mathbb R^d$.  The $m$ parabolic
  equations with constant diffusion coefficient constituting the RDS
  on $\Ot$ are thus pulled-back into equations on a fixed domain,
  albeit with space-time dependent coefficients.  The fixed domain
  setting permits us to use the standard Bochner space machinery
  needed for evolution equations of parabolic type.  On the other
  hand, we are thus left to deal (computationally and analytically)
  with three interacting difficulties: (1) a system of $m$ coupled equations, (2)
  the nature of the nonlinearity $\vec f$ coupling the equations and (3)
  the non-constant diffusion and velocity coefficients, especially as
  functions of time. Our approach in tackling the nonlinearity
  consists in constructing a suitable globally Lipschitz extension of
  the nonlinear reaction field that coincides with it in a
  neighbourhood of the exact solution and then proving that both the
  exact solution and the numerical solution are confined to the domain
  of the original (non-extended) nonlinearity.  We use mainly
  parabolic energy techniques, but must have some pointwise control in
  order to bound the nonlinearities.

\additions{
Our treatment of the nonlinear reaction functions is based on the approach  of \cite{thomee1975galerkin}, see also \cite{Crouziex1987,  schatz1980maximum}. An alternative approach to ours would be to construct schemes where an invariant region for the continuous solution \cite{chueh1977positively} is preserved under discretisation, for work in this direction we refer to  \cite{EllStu93, mckenna2007gnn, Hoff1978}.
  }

  Although all our error estimates are derived for the Lagrangian
  formulation, given that the domain evolution is prescribed, they
  carry in a straightforward manner to the Eulerian framework.  The
  situation would be more delicate if the domain evolution was itself
  an unknown as a geometric motion coupled with the RDS, but this is
  outside the scope of this study.
  
  \additions{
  The smooth prescribed evolution case we deal with in this study 
  is of relevance in many applications (see for example \cite{murray2003mathematical}), including but not limited to skin 
  pigment pattern formation during development. We note that in many 
  important applications such as {\it morphogen controlled growth}  where the evolution of 
  the domain is governed by the solution to the RDS \cite{baker2007mechanism},  
  or cell motility \cite{venk11chemotaxis} and tumour growth  \cite{raquelsfem} where the 
  deformation of the cell membrane is governed by a geometric evolution law, the 
  domain itself is an unknown which must be approximated. This more challenging 
  setting warrants further investigation. The transformation to the reference domain, 
  which we make use of in our Lagrangian analysis, would now depend on the solution
  of the RDSs and/or on the geometric properties of the domain leading to the consideration 
  of quasilinear or fully nonlinear RDSs on fixed domains.
 }
  
\subsection{Implicit-explicit schemes}
  The fully discrete method that we analyse, is a fully practical
  me\-thod, implemented in the ALBERTA toolbox (code available upon
  request), using an implicit-explicit backward Euler scheme to derive
  the time-discretisation \citep{ano2006}.  
  
  \additions{
  On fixed domains, Zhang et
  al. \citep{zhang2008second} analyse a second order implicit-explicit
  finite element scheme for the Gray-Scott model and Garvie and
  Trenchea \cite{garvie2007finite} analyse a first order scheme for an
  RDS that models predator prey dynamics.
  In \cite{chaplain2001spatio} the authors propose and briefly analyse a 
  numerical method based on an IMEX time discretisation and spherical harmonics 
  for the spatial approximation of a RDS 
  posed on the surface of a stationary sphere. The a posteriori analysis of finite element 
  methods for RDSs is  treated  in, for example, the book 
   \cite{estep2000estimating}  where systems of coupled parabolic differential equations 
  (reaction-diffusion) and ordinary differential equations are considered and \cite{moore1994posteriori} 
  where 1d scalar quasilinear RDSs are considered. An adaptive finite 
  element method for semilinear RDSs on evolving domains and surfaces 
  is presented in \cite{venkataraman2011adaptive}. Another approach  is   the 
  {\em moving finite element} method, where nodal movement is regarded as an unknown 
  (even on fixed domain problems) and  at each timestep nodes are moved, usually with the 
  goal of controlling the error \cite{miller1981movingII, miller1981moving,baines1994mfe}, for the 
  analysis of the moving finite element method we refer to \cite{dupont1982mesh}. In \cite{Zeg} the authors
  describe an adaptive moving mesh FEM to approximate solutions of  the Gray-Scott RDS 
  on a {\em fixed} domain. 
Recently Mackenzie and
  Madzvamuse \citep{mackenzie2009analysis} analysed a finite
  difference scheme approximating the solution of a linear RDS on a
  domain with continuous spatially linear isotropic evolution.
  }
  
   Our study is novel, in that, we propose and analyse a finite element
  method to approximate RDSs on a domain with continuous (possibly
  nonlinear) evolution.  This creates space-time-dependent
  coefficients impacting the diffusion and the time-derivative term
  which complicates the fully discrete scheme's analysis and requires
  a careful treatment of the timestep, depending on the rate of domain
  evolution.  In spite of it being only first order in time, the
  proposed implicit-explicit method is robust for the applications we
  have in mind, where long time integration is essential and
  the problems are often posed on complex geometries such as the
  surface of an organism.
\subsection{Outline}
\additions{
The structure of this paper is as follows: in \S\ref{secn.:setup} we
introduce the notation employed throughout this article, we state our
model problem together with the assumptions that we make on the
problem data and the domain evolution. We present the weak formulation
of the continuous problem and define a modified nonlinear reaction
function which we introduce for the analysis. In \S\ref{secn.:fem} we
present the \emph{semidiscrete} (space-discrete) and the \emph{fully
  discrete finite element schemes} with some remarks regarding
implementation, allowing the practically minded reader to skip over
the analysis through to \S\ref{secn.:implementation}.  We then analyse
the semidiscrete scheme in \S\ref{sec:analys-semid-scheme} and the
fully discrete scheme in \S\ref{subsecn.:fully_disc} proving optimal
rate error bounds as well as a maximum-norm stability result, whereby
the stabilising effect of domain growth observed in the continuous
case is preserved at the discrete level and in the numerical schemes.
In \S\ref{secn.:implementation} we provide a concrete
implementation of the finite element scheme
with a set of reaction kinetics commonly encountered in developmental
biology, considering domains with spatially linear and nonlinear
evolution. In \S\ref{secn.:numerics} we present computational
experiments to illustrate our theoretical results.
}
\section{Notation and Setup}\label{secn.:setup}
In this section we define most of the basic notation for the rest of
the paper, introduce the evolving domain framework, set the detailed
blanket assumptions and introduce a pulled-back version of Problem
\ref{Pbm:RDS_td}.
\subsection{Calculus and function spaces}
Given an open and bounded stationary domain $\Pi\subset\Reals^{d}$ and a function $\vec{\eta}\in{C}^{1}(\Pi;\Reals^m)$, we denote by $\nabla\vec\eta$ the Jacobian matrix of $\vec{\eta}$ with components $\left[\nabla\vec\eta(\vec x)\right]_{i j}=\partial_{x_i}\eta_j$. For $\vec{\eta}\in{C}^{1}(\Pi;\Reals^d)$  we denote by $\nabla\cdot\vec\eta$ the divergence of $\vec{\eta}$.
In an effort to compress notation for spatial derivatives, we introduce the convention used above, that if the variable with respect to which we differentiate  is omitted, it should be understood as the spatial argument of the function.
\condense{\begin{equation}
\nabla\vec{\eta}(\vec{x}):=\left[\begin{array}{ccc}
\partial_{x_1}\eta_1(\vec{x})&\dotsc&\partial_{x_1}\eta_m(\vec{x})\\					\vdots		     &\ddots&		\vdots  \\
\partial_{x_d}\eta_1(\vec{x})&\dotsc&\partial_{x_d}\eta_m(\vec{x})\end{array}\right], \text{ for } \vec{x}\in\Pi, 
\end{equation}
and the divergence of $\vec{\eta}$
\begin{align}
\nabla\cdot\vec{\eta}(\vec{x}):=\sum_{i=1}^d\partial_{x_i}{\eta}_i(\vec{x}).
\end{align}
}
\condense{
For the case of a scalar-valued function ${\eta}\in{C}^{1}(\Pi;\mathbb{R})$, we define the Laplacian of $\eta$
\begin{align}
\lap\eta(\vec{x}):=\sum_{i=1}^d\partial_{x_ix_i}\eta(\vec{x}).
\end{align}
}

We denote by $\Lp{p}{(\Pi)}$,  $\Sob{p}{k}{(\Pi)}$ and $\Hil{k}{(\Pi)}$ the Lebesgue, Sobolev and Hilbert spaces respectively,\changesvtwo{equipped with the usual norms and seminorms \citep{evans2009partial}.}For vector valued functions $\vec{\eta},\vec{\mu} : \Pi\to\Reals^m$, we denote
 \begin{align}
\ltwop{\vec{\eta}}{\vec{\mu}}{\Pi^m}&:=\sum_{i=1}^m\int_{\Pi}\eta_i(\vec{x})\mu_i(\vec{x})\dif x,
\end{align}
with the corresponding modifications to the norms and seminorms.

\condense{
as defined by the following: for measurable $\eta$ and for $p,k\in[1,\infty)$
\begin{align}
\Lp{p}{(\Pi)}&:=\left\{\eta:\int_{\Pi}\lv\eta\rv^p<\infty\right\},\\
\Lp{\infty}{(\Pi)}&:=\left\{\eta:\sup_{\vec{x}\in\Pi}\lv\eta(\vec{x})\rv<\infty\right\},\\
\Sob{k}{p}{(\Pi)}&:=\left\{\eta\in\Lp{p}{(\Pi)}:\sum_{\alpha\leq k}\nabla^\alpha\eta\in\Lp{p}{(\Pi)}\right\},\\
\Hil{k}{(\Pi)}&:=\Sob{2}{k}{(\Pi)}.
\end{align}
For measurable functions $\eta,\mu : \Pi\to\Reals$,  we introduce the following notation 
\begin{align}
\ltwop{\eta}{\mu}{\Pi}&:=\int_{\Pi}\eta(\vec{x})\mu(\vec{x})\dif x,\\
\ltwon{\eta}{(\Pi)}&:=\ltwop{\eta}{\eta}{\Pi}^{1/2},\\
\seminorm{\eta}{\Hil{k}(\Pi)}&:=\ltwon{\nabla^{k}\eta}{(\Pi)}, \text{ for } k\in\Integers_{+},\\
\Hiln{\eta}{k}{(\Pi)}&:=\left(\ltwon{\eta}{(\Pi)}^2+\sum_{j=1}^k\seminorm{\eta}{\Hil{j}(\Pi)}^2\right)^{1/2}.
\end{align}}

\subsection{Evolving domain}
\label{sec:evolving-domain}
\additions{
Let $\Oc\subset\mathbb R^d$ be a simply connected, convex domain with Lipschitz boundary;
we will call it the \emph{reference domain}.  We define the
\emph{evolving domain} as a time-parametrised family of domains
}
\begin{equation}
  \label{eqn:mapping_def}
  \{\Ot:=\vec\A_t(\Oc)\}_{0\leq t\leq T}
  \text{ where }
  \vec\A_t:\Oc \rightarrow\Omega_t
  \text{ is a $\Cont1$-diffeomorphism for each fixed }
  t\in[0,T].
\end{equation}
The Jacobian matrix of $\vec\A_t(\cdot)$, its determinant and inverse will be respectively denoted by
\begin{equation}
  \label{eqn:jacobian_def}
  \vec{J}_t(\vec{\xi})
  :=
  \nabla\A_t(\vec{\xi}) 
  ,
  \quad
  J_t(\vec{\xi}):=\det\vec{J}_t(\vec{\xi})
  \text{ and }
  \vec{K}_t (\vec{\xi}):= [\nabla\A_t(\vec{\xi})]^{-1}
\end{equation}
for each $(\vec \xi,t)\in\Oc\times[0,T]$.
We will use also the \emph{evolution induced convection} on the evolving domain 
\begin{equation}
  \label{eqn:dta}
  \vec a(\vec x,t)
  :=
  \partial_t\vec{\A}_t(\vec \A_t^{-1}(\vec x))
  \text{ for }
  \vec x\in\O_t\mbox{ and } t\in[0,T].
\end{equation}
From classical results \citep{acheson1990elementary} we have the following expression 
\begin{equation}\label{eqn:dtdetj}
  \pdt{J}\left(\vec{\xi},t\right)=J_t(\vec{\xi})\nabla\cdot\vec{a}\left(\A_t(\vec{\xi}),t\right)
  \text{ for } 
  (\vec \xi, t)\in\Oc\times[0,T],
\end{equation} 
\additions{
and the Reynold's transport theorem \citep{acheson1990elementary}, which reads: For a function $g\in C^1(\Ot,[0,T])$
\begin{equation}\label{eqn:transport-theorem}
\frac{\dif}{\dif t}\int_{\Ot}g=\int_{\Ot}\pdt{g}+\nabla\cdot\left(\vec{a}g\right).
\end{equation}
}
To aid the exposition we define $\Q$ to be the topologically cylindrical  space-time domain: 
\begin{equation}
\Q:=\left\{(\vec{x},t):\vec{x}\in\Ot,t\in[0,T]\right\}.
\end{equation}
We now introduce notation to relate functions defined on the evolving domain to functions defined on the reference domain. Given a function $g:\Q\rightarrow\mathbb{R}$ we denote by $\hat{g}:\Oc\times[0,T]\rightarrow\mathbb{R}$  its pullback on the reference domain, defined by the following relationship
\begin{align}\label{eqn:evolving_eulerian_correspondence}
\hat{g}(\vec{\xi},t):=g\left(\A_t(\vec{\xi}),t\right) \quad (\vec \xi, t)\in\Oc\times[0,T].
\end{align}
Assuming sufficient smoothness on the function $g$, using
(\ref{eqn:evolving_eulerian_correspondence}) and the chain rule we may
relate time-differentiation on the reference and evolving
domains:\footnote{To avoid confusion, as in
  \eqref{eqn:material_derivative}), we denote by
  $\partial_if$ the partial derivative with respect to the
  $i$-th argument of the function $f$, for a positive integer $i$.  
  When there is no risk of
  confusion we write $\partial_t f$ for the time derivative of a
  time-dependent function $f$ even when such a variable is not
  explicitly written in the arguments.}
\begin{align}\label{eqn:material_derivative}
\pdt\hat{g}(\vec{\xi},t)=&
\partial_2{g} \left(\A_t(\vec{\xi}),t\right) +\left[\vec{a}\cdot\nabla{g}\right]\left(\A_t(\vec{\xi}),t\right), \quad (\vec \xi, t)\in\Oc\times[0,T].
\end{align} 
The right hand
side of (\ref{eqn:material_derivative}) is commonly known as the
material derivative of $g$ with respect to the velocity $\vec{a}$. The
following result relates the norm of a function
$g:\Q\rightarrow\Reals$ on the evolving domain with its pullback $\hat
g$ on the reference domain:
\begin{align}\label{eqn:norm_evolving_eulerian_correspondance}
\ltwon{g}{(\Ot)}^2&=\ltwop{J_t\hat{g}}{\hat{g}}{\Oc}=:\mynorm{\hat{g}}{J_t}^2.
\end{align}
For the gradient of a sufficiently smooth function $g:\Q\rightarrow\Reals$, we have
\begin{align}\label{eqn:gradient_evolving_eulerian_correspondence}
  \ltwon{\nabla{g}}{(\Ot)}^2&=\ltwop{J_t\vec K_t\nabla\hat{g}}{\vec K_t\nabla\hat{g}}{\Oc}
  =\ltwop{\vec B_t\nabla\hat{g}}{\nabla\hat{g}}{\Oc}
  =:\seminorm{\hat g(\cdot,t)}{\vec B_t}^2,
\end{align}
where $\vec B:=J\vec K\vec K^\transpose$.
For $t\in[0,T]$ we define the bilinear form 
\begin{equation}
  \bbil{\hat{v}}{\hat{w}}{\Oc}
  :=
  \langle \vec{B}_t\nabla \hat{v},\nabla \hat{w} \rangle_{\Oc}
  \text{ for }\hat{v},\hat{w}\in\Hil1(\Oc).
\end{equation}
Assumption \ref{assum:mapping} implies that there exists $\mu,\bar{\mu}\in\Reals_+$ such that for $i=1,\dotsc,m$ and for all $\hat{v}\in\Hil{1}(\Oc)$,
\begin{equation}\label{eqn:coercivity}
\mu\ltwon{\nabla\hat{v}}{(\Oc)}^2\leq \bbil{\hat v}{\hat v}{\Oc} \leq \linfn{\vec{B}}{(\Oc)}\ltwon{\nabla\hat v}{(\Oc)}^2= \bar{\mu}\ltwon{\nabla\hat v}{(\Oc)}^2.
\end{equation}
\begin{Assum}[Regularity of the mapping]\label{assum:mapping}
  \additions{
    It will be handy sometimes to denote the family
    $\setsi{\vec{\A}}t{\clinter0T}$, introduced in
    \ref{eqn:mapping_def}, $\vec\A$ as a single map
    $(\vec{x},t)\ni\W\times[0,T]\mapsto\vec{\A}_t(\vec{x})$.  We assume
    the following regularity:
  }
  \begin{equation}
    \A\in\Cont{1}(\Oc\times[0,T])\quad\text{ and }\quad\A_t\in\Cont{k+1}(\Oc)\text{ for each }t\in[0,T],
  \end{equation}
  where $k$ will be taken equal to the degree of the basis functions of
  the finite element space defined in the following section. To ensure
  the mapping is invertible we assume the determinant of the Jacobian
  $J$ of the mapping $\A$ (cf. (\ref{eqn:jacobian_def})) satisfies
  \begin{equation}
    J>0\mbox{ in }\Oc\times[0,T].\label{eqn:positive_jacobian}
  \end{equation}
\end{Assum}
\subsection{The RDS reformulated on the reference domain}
\label{Pbm:RDS_rd}
Using
(\ref{eqn:evolving_eulerian_correspondence})---(\ref{eqn:gradient_evolving_eulerian_correspondence})
and a change of variable in the divergence, we obtain the following
equivalent formulation of Problem \ref{Pbm:RDS_td} on a reference
domain. Denote by $\hat{\vec u}:\Oc\times[0,T]\to\Reals^m$ the function 
that satisfies \fim
\begin{equation}\label{eqn:reference_model_problem}
   \begin{gathered}
     \pdt\hat{u}_i
     -\frac{D_i}{J}\nabla\cdot(
     \vec{B}
     \nabla\hat{u}_i
     )
     +\hat{u}_i
     \nabla\cdot\hat{\vec a}
     =f_{i}\left(\hat{\vec{u}}
     \right), 
     \text{ on }
     \Oc\times(0,T],
     \\
     \hat{\normal}\cdot\vec{B}\nabla{\hat{u}}_i=0, 
     \text{ on }
     \partial\Oc\times(0,T],
     \\
     \hat{u}_i(\vec{\xi},0)=\hat{u}_i^{0}(\vec{\xi}), 
     \quad
     \vec{\xi}\in\Oc.
   \end{gathered}
\end{equation}
where
\begin{equation}
  \hat{\vec a}(\A_t(\vec\xi))=\partial_t\A_t(\vec\xi)
  \text{ for }
  (\vec\xi,t)\in\Oc\times[0,T].
\end{equation}
 \additions{ 
\begin{Assum}[Nonlinear reaction vector field]
  \label{assum:reaction}
  We assume throughout that $\vec f$ is of the form
  \begin{equation}\label{eqn:f_structure}
    {f}_i(\vec{z})={z}_i F_i(\vec{z}),
    \text{ for all }\vec z\in\operatorname{Dom}\vec f=:I
    \text{ and each }i=1,\dotsc,m
    , 
  \end{equation}
  for some vector field $\vec F\in\Cont1(I)$ and some open set $I\subset\Reals^m$.  As a result $\vec f\in\Cont1(I)$ and 
  it is locally Lipschitz. 
In \S \ref{secn.:implementation} we provide an example of a widely studied set of reaction kinetics that satisfy the structural assumptions we make on the nonlinear reaction vector field. 
\end{Assum}
  }
\begin{Assum}[Existence and regularity]
  \label{assum:well_posed} 
  We assume the global existence of a solution $\hat{\vec{u}}$ to
  Problem \ref{Pbm:RDS_rd}. Furthermore we assume $\hat{\vec{u}}$ is
  in $H^{\ell+1}(\Oc)^m$ with $\partial_{t}\hat{\vec{u}}$ in
  $H^{\ell+1}(\Oc)^m$ where $\ell$ is the polynomial degree of the
  finite element space defined in the following section.
\end{Assum}
\begin{Rem}[Applicability of Assumption \ref{assum:well_posed}]
  In \citep{venkataraman2010global}, we proved the global existence of
  positive classical solutions to Problem \ref{Pbm:RDS_td} for a class
  of RDSs with positive initial data on domains with bounded spatially
  linear isotropic evolution. 
  \deletions{
  In this case, if $\vec{f}$ belongs to
  $\Cont1(\Reals^m_+)$, the solution $\vec u$ is shown to be contained
  in a rectangle $\prod_{i=1}^d[u_i^-,u_i^+]$, with $\vec u^->0$ and
  assumption \ref{assum:reaction} holds with the region
  ${I}_\delta:=\prod_{i=1}^d[u_i^--\delta,u_i^++\delta]$ a subset of
  $\Reals_+^m$.
  }
\end{Rem}
\subsection{Weak formulation}
To construct a finite element discretisation, we
introduce a {\em weak solution} of Problem
\ref{Pbm:RDS_rd}, denoted by $\hat{u}_i \in
\Lp{2}{\big(0,T;\Hil{1}{(\Oc)}\big)}$ with
$\pdt\hat{u}_i\in\Lp{2}{\big(0,T;\Hil{-1}{(\Oc)}\big)}$ such that
\begin{equation}
\label{eqn:weak_form_reference_jddtu_uddtj}
\\
\ltwop{
    J\qp{\pdt\hat{u}_i
      +
      \hat u_i\nabla\cdot\vec{a}
      \qp{
        \A_t(\cdot,t)
      }
  }}{
    \hat{\chi}
  }{\Oc}
  +
  D_i
  \bbil{
    {\nabla}\hat{u}_i
  }{
    \nabla\hat{\chi}
  }{\Oc}
  =
  \ltwop{J{f}_i(\hat{\vec{u}})}{\hat{\chi}}{\Oc} 
  \Foreach\hat{\chi} \in \Hil{1}{(\Oc)}.
\end{equation}
Using the expression for the time-derivative of the determinant of the Jacobian (\ref{eqn:dtdetj}), we have 
\begin{equation}
\begin{split}
\label{eqn:cont_weak_form_reference_ddtju}
\ltwop{\pdt(J\hat{u}_i)}{\hat{\chi}}{\Oc}+{D}_i\bbil{{\nabla}\hat{u}_i}{\nabla\hat{\chi}}{\Oc}&=\ltwop{J f_i(\hat{\vec{u}})}{\hat{\chi}}{\Oc}\Foreach\hat{\chi} \in \Hil{1}{(\Oc)}.
\end{split}
\end{equation}
We shall use (\ref{eqn:cont_weak_form_reference_ddtju}) to construct a finite element scheme to approximate the solution to Problem \ref{Pbm:RDS_td} on the reference domain. 
\subsection{Extended nonlinear reaction function}
\additions{In general the techniques used to show Assumptions
  \ref{assum:reaction} and \ref{assum:well_posed} hold utilise the
  maximum principle \citep{smoller1994swa, venkataraman2010global}. In
  the discrete case, since the maximum principle cannot be applied
  \citep[p. 83]{Thomee} we show, under suitable assumptions,
  maximum-norm bounds on the discrete solution in
  (\ref{eqn:linf_closeness}) that guarantee the solution remains in
  the region ${I}$ defined in \ref{eqn:f_structure}.  We introduce a
  modified {\em globally Lipschitz} nonlinear reaction in order to
  derive the error bounds, but this extension is never needed in
  practice and hence needs not be computed.  Recalling
  Assumption \ref{assum:reaction}, we define
  $\widetilde{\vec{F}}\in{C}^1(\mathbb{R}^m)$ such that
\begin{equation}\label{eqn:ftilde_def}
  \begin{gathered}
    \begin{cases}
      \begin{aligned} 
        \widetilde{\vec{F}}(\vec{z})&=\vec{F}(\vec{z}), \quad \mbox{for } \vec{z}\in{I}\\
        \lv\widetilde{\vec{F}}^\prime(\vec{z})\rv &< \widetilde{C}, \quad  \mbox{for } \vec{z}\in\mathbb{R}^m,
      \end{aligned}
    \end{cases}
    \qquad
    \text{ and }
    \widetilde{{f}_i}(\vec{z})
    :=
    z_i{\widetilde{F}_i}(\vec{z}), \quad \mbox{for } \vec{z}\in\mathbb{R}^m.
  \end{gathered}
\end{equation}
The function $\widetilde{\vec F}$  is guaranteed to exist due to Assumptions
\ref{assum:reaction}, \ref{assum:well_posed}, the Whitney Extension
Theorem \citep[Th. 1, \S 6.5]{EvansGariepy:92} and the use of an
appropriate cut-off factor.  If $\vec{u}$ is a solution of
(\ref{eqn:model_problem})
\begin{equation}\label{eqn:f_tilde_equals_f_for_exact_solution}
  \begin{split} 
    \widetilde{\vec{f}}(\vec{u})&={\vec{f}}(\vec{u}). 
  \end{split}
\end{equation}
Thus, we may without restriction replace $\vec{f}$ with $\widetilde{\vec{f}}$ in (\ref{eqn:model_problem})}
\section{Finite element method}\label{secn.:fem}
In this section we design the finite
element method, first by discretising Problem \ref{Pbm:RDS_rd} in
space only, discussing some properties of the semidiscrete scheme and
then passing to the fully discrete scheme.
\subsection{Spatial discretisation set-up}
We shall split the spatial and temporal discretisation of Problem
\ref{Pbm:RDS_td} into separate steps. For the spatial approximation,
we employ a conforming finite element method.  To this end, we define
$\Th$ a triangulation of the reference domain. We shall consistently
denote by $\hat{h}:=\max_{\s\in\Th}\operatorname{diam}(\s)$ the
mesh-size of $\Tc$.\changesvtwo{We assume the triangulation $\Th$ is conforming and that there is no error due to boundary approximation.
Furthermore given $\lbrace\Th_i\rbrace_{i=1}^\infty$, a sequence of conforming
triangulations, we assume the {\em quasi-uniformity} of the sequence
holds, for details see for example  \citep{schwab1998p}. Note the assumption of
quasi-uniformity  implies that the family of triangulations is
{shape-regular} \citep[p. 159]{schwab1998p}.}
\condense{We assume the triangulation $\Th$ fulfils the
following properties:
\begin{itemize}
\item By $\s \in \Th$ we mean $\s$ is an open simplex (interval,
  triangle or tetrahedron for $d=1,2$ or $3$ respectively).
\item $\Th$ is conforming, i.e., for any $\s,k\in\Th$,
  $\bar{\s}\cap\bar{k}$ is either $\¤emptyset$, a vertex, an edge or a
  face common to $\s$ and $k$ or the full simplex $\bar{\s}=\bar{k}$.
\item No error due to boundary approximation, i.e.,
  $\cup_{\s\in\Th}\bar{\s}=\bar{\Oc}$ (we make this assumption for
  ease of exposition and it may be relaxed depending on the
  application).
\end{itemize} 
For a sequence $\lbrace\Th_i\rbrace_{i=1}^\infty$ of conforming
triangulations, we assume the {\em quasi-uniformity} of the sequence
holds, i.e., there exist $C_1,C_2$ independent of $i$ such that
\begin{equation}
  C_1\hat{h}\leq\hat{h}_s\leq{C}_2\hat{h}_\s^\circ, \quad \text{ for all } \s\in\Th_i, i=1,2,\dotsc,
\end{equation}
where $\hat{h}_\s^\circ$ and $\hat{h}_\s$ are the radius of the
largest ball contained in $\s$ and the diameter of $\s$
respectively. Furthermore we note that our assumption of
quasi-uniformity implies that the family of triangulations is
{shape-regular} \citep[p. 159]{schwab1998p}.
}

Given the triangulation $\Th$, we now define a {finite element space}
on the reference configuration:
\begin{equation}
  \label{eqn:fe_space_ref}
  \Vc := \left\{ \hat{\Phi} \in H^1(\Oc) :  \hat{\Phi}|_s \mbox{ is piecewise polynomial of degree $\ell$}\right\}. 
\end{equation}
We utilise the following known results about the accuracy of the
finite element space $\Vc$. By the definition of $\Vc$, we have for
$\hat{v} \in H^{\ell+1}(\Oc)$ (see for example
\citet{brenner2002mathematical} or \citet{Thomee}),
\begin{equation}
  \begin{split}\label{eqn:reference_fe_space_accuracy}
    \inf_{\hat{\Phi}\in\Vc}\left\{\ltwon{\hat{v}-\hat{\Phi}}{(\Oc)}+\hat{h}\ltwon{\nabla(\hat{v}-\hat{\Phi})}{(\Oc)}\right\} \leq C\hat{h}^{\ell+1}\lv\hat{v}\rv_{H^{\ell+1}(\Oc)}.
  \end{split}
\end{equation}
\additions{
Let  the degree of the finite element space satisfy $\ell+1 >
\frac{d}{2}$, where $d$ is the spatial dimension.
In the analysis we shall make use of the fact that
(\ref{eqn:reference_fe_space_accuracy}) is satisfied by taking the
Lagrange interpolator $\Lagrange:\Hil{\ell+1}{(\Oc)}\to\Vc$ in
place of $\hat{\Phi}$ (note that  $\ell+1 >
{d}/{2}$ implies $\Hil{\ell+1}(\Oc)\hookrightarrow \Cont{0}(\Omega)$ so the Lagrange interpolant is well defined). Let $\mathcal{I}^h:\Cont0\rightarrow\Vc$ be a
Cl\'ement type interpolant \citep{Clement:1975}. The following bound
holds
}
\begin{equation}
  \begin{split}\label{eqn:reference_fe_space_linf_accuracy}
    \|\hat{v}-\mathcal{I}^h\hat{v}\|_{\Lp{\infty}(\Oc)}\leq{C}\hat{h}^{\ell+1-d/2}\lv\hat{v}\rv_{H^{\ell+1}(\Oc)}.
  \end{split}
\end{equation}
We shall make use of the following inverse estimate valid on
quasiuniform sequences of triangulations:
\begin{equation}
  \begin{split}
    \label{eqn:fe_space_reference_inv_est}\|\hat{\Phi}\|_{\Lp{\infty}(\Vc)} \leq C\hat{h}^{-d/2}\|\hat{\Phi}\|_{\Lp{2}(\Vc)} \quad \forall \hat{\Phi} \in \Vc.
  \end{split}
\end{equation}
\subsection{Semidiscrete approximation}\label{subsecn.:semi_disc}
We define the spatially semidiscrete approximation of the solution of
Problem \ref{Pbm:RDS_td} to be a function
$\hat{u}^h_i:[0,T]\rightarrow\Vc$, such that for $i=1,\dotsc,m$,
\begin{equation}
  \begin{cases}
    \begin{aligned}
      \label{eqn:sd_weak_form}
      \ltwop{\pdt(J\hat{u}^h_i)}{\hat{\Phi}}{\Oc}+\ltwop{{D}_i\vec{B}{\nabla}\hat{u}^h_i}{\nabla\hat{\Phi}}{\Oc}&=\ltwop{J\widetilde{f}_i(\hat{\vec{u}}^h)}{\hat{\Phi}}{\Oc}
      \Foreach\hat{\Phi} \in \Vc,\\
      \hat{u}^h_i(0)&=\Lambda^h\hat{u}_i^0, 
    \end{aligned}
  \end{cases}
\end{equation}
where $\Lagrange$ is the Lagrange interpolant.
\begin{Prop}[Solvability of the semidiscrete scheme]\label{Lem:semidiscrete_solvability}
  Let Assumptions \ref{assum:well_posed} and \ref{assum:mapping}
  hold. Then, the semidiscrete scheme (\ref{eqn:sd_weak_form})
  possesses a unique solution
  $\hat{\vec{u}}^h\in\Lp{\infty}{(0,T)}^m$.
\end{Prop}
\begin{Proof}
  In (\ref{eqn:sd_weak_form}) if we write $\hat{u}^h_i(t)$ as
  $\sum_{j=1}^{\dim(\Vc)}\alpha_j\hat{\Phi}_j$, we obtain a system of
  $\dim(\Vc)$ ordinary differential equations for each $i$. By
  assumption the initial data for each ODE is bounded. From Assumption
  \ref{assum:mapping} and the construction of $\widetilde{\vec{f}}$
  (\ref{eqn:ftilde_def}), we have that $J$, $\widetilde{\vec{f}}$ and
  their product are continuous globally Lipschitz functions. From ODE
  theory (for example \citep{schmitt1998nonlinear}) we conclude that
  (\ref{eqn:sd_weak_form}) possesses a unique bounded solution.
\end{Proof}
\subsection{The effect of domain evolution on the semidiscrete solution}\label{subsecn.:stabilising_effect_of_growth}
We now examine the stability of (\ref{eqn:sd_weak_form}) and show that
domain growth has a {\em diluting} or {\em stabilising} effect on the
semidiscrete solution, mirroring results for the continuous problem
\citep{labadie2007stabilizing}. Taking $\hat{\Phi}=\hat{u}^h_i$ in
(\ref{eqn:sd_weak_form}) gives \fim
\begin{equation}
  \label{eqn:stability_weakform_I}
  \ltwop{\pdt(J\hat{u}^h_i)}{\hat{u}^h_i}{\Oc}+{D}_i\bbil{{\nabla}\hat{u}^h_i}{\nabla\hat{u}^h_i}{\Oc}=\ltwop{J\widetilde{f}_i(\hat{\vec{u}}^h)}{\hat{u}^h_i}{\Oc}.
\end{equation}
For the first term on the left of ({\ref{eqn:stability_weakform_I}})
we have
\begin{equation}
\ltwop{\pdt(J\hat{u}^h_i)}{\hat{u}^h_i}{\Oc}=\frac{\dif}{\dif t}\ltwop{J\hat{u}^h_i}{\hat{u}^h_i}{\Oc}-\ltwop{J\hat{u}^h_i}{\pdt\hat{u}^h_i}{\Oc}.
\end{equation}
\additions{Application of Reynold's transport theorem
(\ref{eqn:transport-theorem}) gives
\begin{equation}\label{eqn:pdtJhatu_u_on_evolving_frame}
  \ltwop{\pdt(J\hat{u}^h_i)}{\hat{u}^h_i}{\Oc}
  =
  \frac{1}{2}\left(\frac{\dif}{\dif t}\mynorm{\hat{u}^h_i}{J}^2
  +
  \ltwop{J{u}^h_i}{{u}^h_i\nabla\cdot\vec{a}(\A_t(\vec \xi),t)}{\Oc}\right).
\end{equation}}
Dealing with the right hand side of
(\ref{eqn:stability_weakform_I}), using (\ref{eqn:ftilde_def}) and the
mean-value theorem (MVT) we have with $\widetilde{C}$ from
(\ref{eqn:ftilde_def})
\begin{equation}\label{eqn:ftilde_MVT}
  \begin{split}
    \lv\widetilde{f}_i(\hat{\vec{u}}^h)\rv&\leq\lv\widetilde{f}_i(\vec{0})\rv
    +
    \lv\widetilde{f}_i(\hat{\vec{u}}^h)-\widetilde{f}_i(\vec{0})\rv
    \leq\lv\widetilde{f}_i(\vec{0})\rv
    +
    \widetilde{C}\sum_{j=1}^m\lv\hat{u}_j^h\rv
    .
  \end{split}
\end{equation}
Therefore we have
\begin{equation}
  \lv\ltwop{J\widetilde{f}_i(\hat{\vec{u}}^h)}{\hat{u}^h_i}{\Oc}\rv
  \leq
  \widetilde{C}\ltwop{J\sum_{j=1}^m\lv\hat {u}^h_j\rv}{\lv \hat {u}^h_i\rv}{\Oc}+\lv\ltwop{J \widetilde{f}_i(0)}{\hat u^h_i}{\Oc}\rv.
\end{equation}
Applying Young's inequality gives
\begin{equation}
\begin{split}
\lv\ltwop{J\widetilde{f}_i(\hat{\vec{u}}^h)}{\hat{u}^h_i}{\Oc}\rv
\leq
&\widetilde{C}\left(\frac{1}{2}\sum_{j\not= i}\mynorm{\hat{u}^h_j}{J}^2
 +\frac{m+1}{2}\mynorm{\hat{u}^h_i}{J}^2\right)+\frac{1}{2}\mynorm{\hat{u}^h_i}{J}^2+C_{\widetilde{f}_i(\vec{0})},
\end{split}
\end{equation}
where $C_{\widetilde{f}_i(\vec{0})}\in\Reals^+$ depends on $\lv\widetilde{f}_i(\vec{0})\rv$.  Summing over $i$ we have
\begin{equation}\label{eqn:J_vecftilde_stability_estimate}
\begin{split}
\sum_{i=1}^m\lv\ltwop{J\widetilde{f}_i(\hat{\vec{u}}^h)}{\hat{u}^h_i}{\Oc}\rv&
\leq
\left(\widetilde{C}m+\frac{1}{2}\right)\mynorm{\hat{\vec{u}}^h}{J^m}^2+C_{\widetilde{\vec{f}}(\vec{0})}.
\end{split}
\end{equation}
Using (\ref{eqn:gradient_evolving_eulerian_correspondence}), (\ref{eqn:pdtJhatu_u_on_evolving_frame}) and (\ref{eqn:J_vecftilde_stability_estimate}) in (\ref{eqn:stability_weakform_I}) gives
\begin{equation}
\begin{split}
\frac{\dif}{\dif t}\mynorm{\hat{\vec{u}}^h}{J^m}^2+2\sum_{i=1}^m D_i\seminorm{\hat{u}^h_i}{\vec B}^2
\leq
&\ltwop{J\left(2\widetilde{C}m+1-\nabla\cdot\vec{a}\left(\A_t(\vec{\xi}),t\right)\right)\hat{\vec{u}}^h}{\hat{\vec{u}}^h}{\Oc^m}+2C_{\widetilde{\vec{f}}(\vec{0})}.
\end{split}
\end{equation}
Finally, integrating in time and applying Gronwall's lemma we have
\begin{equation}\label{eqn:growth_estimate}
\mynorm{\hat{\vec{u}}^h(t)}{J^m}^2
\leq
\left(\mynorm{\hat{\vec{u}}^h(0)}{J^m}^2+2t C_{\widetilde{\vec{f}}(\vec{0})}\right)\exp\left(\sup_{\Oc\times[0,T]}\left\{2\widetilde{C}m+1-\nabla\cdot\vec{a}\left(\A_t(\vec{\xi}),t\right)\right\}t\right).
\end{equation}
From (\ref{eqn:dtdetj}), the dilution term
$\nabla\cdot\vec{a}$ has the same sign as $\pdt{J}$ and is therefore
positive (or negative) if the domain is growing (or contracting). Thus,
domain growth has a diluting effect on the $\Lp{2}{(\Ot)^m}$ norm (c.f., (\ref{eqn:norm_evolving_eulerian_correspondance})) of
the solution.
\subsection{Fully discrete scheme}
We divide the time interval $[0,T]$ into $N$ subintervals, $0 = t_0 <
\dotsb < t_N = T$ and denote by $\tau_n := t_{n} - t_{n-1}$ the
(possibly nonuniform) time step and $\tau=\max_n{\tau_n}$. We
consistently use the following shorthand for a function of time: $ f^n
:= f (t_n )$, we denote by
$\bar{\partial}f^n:={\tau_n}^{-1}\left(f^n-f^{n-1}\right).$

For the approximation in time we use a modified implicit Euler method
where linear reaction terms and the diffusive term are treated
implicitly while the nonlinear reaction terms are treated
semi-implicitly using values from the previous timestep (the first step
of a Picard iteration). Our choice of timestepping scheme stems from
the numerical investigation conducted by \citet{ano2006}.

The fully discrete scheme we employ to approximate the solution of
Problem \ref{Pbm:RDS_td} is thus, find $\hat{U}_i^n\in\Vc$, for
$n=1,\dotsc,N$, such that for $i=1\dotsc,m$, we have
\begin{equation}\label{eqn:fully_discrete_scheme}
\begin{cases}
\begin{split}
\ltwop{\bar{\partial}\left[J\hat{U}_i\right]^n}{\hat{\Phi}}{\Oc}+{D}_i&\ltwop{[\vec{B}\nabla\hat{U}_i]^n}{\nabla\hat{\Phi}}{\Oc}=\ltwop{J^n\hat{U}_i^n\widetilde{F}_i(\hat{\vec{U}}^{n-1})}{\hat{\Phi}}{\Oc}, \quad \forall\hat{\Phi}\in\Vc,\\
\hat{U}_i^0&=\Lambda^h\hat{u}_i^0,
\end{split}
\end{cases}
\end{equation}
where $\Lagrange$ is the Lagrange interpolant and $\widetilde F_i$ is as defined in (\ref{eqn:ftilde_def}). 
\deletions{
\subsection{Moving mesh formulation}
Alternatively, and in a more physically intuitive way, we may look to
approximate the solution to (\ref{eqn:model_problem}) on a conforming
subspace of the evolving domain. To this end we define a family of
finite dimensional spaces $\Vtn, n= [0,\dotsc,N]$ such that with
$\hat{\Phi}$ as defined in (\ref{eqn:fe_space_ref}):
\begin{equation}
\label{eqn:fe_space_mov}
\Vtn := \left\{ {\Phi^n}\in\Hil{1}(\Otn):{\Phi^n}(\vec{\A}_{t^n}(\vec{\xi}))=\hat{\Phi}(\vec{\xi})\right\},
\end{equation}
which also defines the \emph{triangulation} $\T^n$, $n= [0,\dotsc,N]$  on the evolving domain. Using (\ref{eqn:fully_discrete_scheme}) and (\ref{eqn:fe_space_mov}) we have the following equivalent finite element formulation on the evolving domain: find ${U}_i^n\in\Vtn$, for $n=1,\dotsc,N$, such that for $i=1\dotsc,m$,
\begin{equation}\label{eqn:moving_fully_discrete_scheme}
\begin{cases}
\begin{aligned}
\bar{\partial}\left[\ltwop{{U}_i}{{\Phi}}{\Ot}\right]^n+{D}_i\ltwop{\nabla{U}_i^n}{\nabla{\Phi}^n}{\Otn}&=&&\ltwop{{U}_i^n\widetilde{F}_i({\vec{U}}^{n-1})}{{\Phi}^n}{\Otn},\\
{U}_i^0&=&&\Lambda^h{u}_i^0,
\end{aligned}
\end{cases}
\end{equation}
where $\Lambda^h:\Hil{1}{(\O_0)}\to\V^0$ is the Lagrange interpolant.
By (\ref{eqn:material_derivative}) and (\ref{eqn:fe_space_mov}) 
\begin{equation}
\pdt{\hat{\Phi}}(\vec{\xi})=\left[{\pdt\Phi}+\vec{a}\cdot\nabla\Phi\right]\left(\vec{\A}_t\left(\vec{\xi}\right)\right)=0.
\end{equation}
}
\additions
{
\subsection{Physical domain formulation}
In a more physically intuitive way, we may look to
approximate the solution to (\ref{eqn:model_problem}) on a conforming
subspace of the evolving domain. To this end we define a family of
finite dimensional spaces $\Vtn, n= [0,\dotsc,N]$ such that 
\begin{equation}
\label{eqn:fe_space_mov}
\Vtn := \left\{ \hat\Phi(\A_{t^n}^{-1}(\cdot)):\hat\Phi\in\Vc\right\},
\end{equation}
which also defines the \emph{triangulation} $\T^n$, $n= [0,\dotsc,N]$  on the evolving domain. Using (\ref{eqn:fully_discrete_scheme}) and (\ref{eqn:fe_space_mov}) we have the following equivalent finite element formulation on the evolving domain: Find ${U}_i^n\in\Vtn$, for $n=1,\dotsc,N$, such that for $i=1\dotsc,m$,
\begin{equation}\label{eqn:moving_fully_discrete_scheme}
\begin{cases}
\begin{aligned}
\bar{\partial}\left[\ltwop{{U}_i}{{\Phi}}{\Ot}\right]^n+{D}_i\ltwop{\nabla{U}_i^n}{\nabla{\Phi}^n}{\Otn}&=&&\ltwop{{U}_i^n\widetilde{F}_i({\vec{U}}^{n-1})}{{\Phi}^n}{\Otn}\Foreach\Phi^n\in\Vtn\\
{U}_i^0&=&&\Lambda^h{u}_i^0,
\end{aligned}
\end{cases}
\end{equation}
where $\Lambda^h$ is the Lagrange interpolant.
}
\section{Analysis of the semidiscrete scheme}
\label{sec:analys-semid-scheme}
We now prove that the semidiscrete solution converges to the exact one
with optimal order in the
$\Lp{\infty}{(0,T;\Lp{2}{(\Oc)^m})}$ norm
and the $\Lp{2}{(0,T;\Hil{1}{(\Oc)^m})}$
seminorm.
\subsection{A time-dependent Ritz projection}
A central role in the analysis is played by the Ritz, or elliptic
projector, defined, as in \citet{wheeler1973priori}, for each
$t\in[0,T]$, by $\Ritz: H^1(\Oc)\to\Vc$ such that for each
${\hat{v}}\in\Hil{1}{(\Oc)}$
\begin{gather}
  \label{eqn:ritz_reference_definition}
  \bbil{\hat{v} }{\hat{\Phi}}{\Oc}=\bbil{\Ritz\hat{v}}{\hat{\Phi}}{\Oc}\Foreach\hat{\Phi} \in \Vc,
  \\\text{ and }
  \label{eqn:Ritz_constraint_2}
  \int_{\Oc}\left[\Ritz\hat{v}-\hat{v}\right]=0.
\end{gather}
The constraint  (\ref{eqn:Ritz_constraint_2}) ensures $\Ritz$ is well defined. Differentiation in time in (\ref{eqn:ritz_reference_definition}) with $v=\hat{u}_i$ yields
\begin{equation}
  \label{eqn:ritz_dt_defn}
  \bbil{\pdt(\hat{u}_i-\Ritz\hat{u}_i )}{\hat{\Phi}}{\Oc}+\langle (\pdt\vec B)\nabla(\hat{u}_i-\Ritz\hat{u}_i ),\nabla\hat{\Phi}\rangle_{\Oc}=0\Foreach\hat{\Phi} \in \Vc.
\end{equation}
  To obtain optimal error estimates, we now decompose the error into an
{\em elliptic error} (the error between the Ritz projection and the
exact solution) and a {\em parabolic error} (the error between the
semidiscrete solution and the Ritz projection):
\begin{equation}
\begin{split}
  \label{eqn:error_decomp_semi_disc}
  \hat{\vec u}^h-\hat{\vec u}
  &
  =
  (\hat{\vec u}^h-\Ritz\hat{\vec u})
  +
  (\Ritz\hat{\vec u}-\hat{\vec u})
  =:
  \hat{\vec \rho}^h+\hat{\vec \varepsilon}^h.
\end{split}
\end{equation}
where the equality defines $\hat{\vec{\rho}}^h=(\hat{\rho}^h_1,\dotsc,\hat{\rho}^h_m)^\transpose$ and $\hat{\vec{\varepsilon}}=(\hat{\varepsilon}_1,\dotsc,\hat{\varepsilon}_m)^\transpose$.  
\begin{Lem}[Ritz projection error estimate]\label{lem:epsilon_grad_epsilon}
  Suppose assumptions \ref{assum:well_posed} and \ref{assum:mapping}
  (with $k=\ell$) hold and let $\Ritz[]$ be the Ritz projection
  defined in (\ref{eqn:ritz_reference_definition}).  Then the
  following estimates hold:
  \begin{gather}
    \begin{split}
      \label{eqn:Ritz_projection_error}
      \sup_{t\in[0,T]}\Bigg\{&\ltwon{\Ritz\hat{\vec{u}}(t)-\hat{\vec{u}}(t)}{(\Oc)^m}^2
      + \hat{h}^2\sum_{i=1}^m\ltwon{\nabla\left(\Ritz\hat{{u}}_i(t)-\hat{{u}}_i(t)\right)}{(\Oc)}^2\Bigg\}
\leq
{C(\A,\hat{u})}\hat{h}^{2(\ell+1)},
    \end{split}
    \\
    \label{eqn:dt_Ritz_projection_error}
    \begin{split}
      \sup_{t\in[0,T]}\Bigg\{&\ltwon{\pdt\left(\Ritz\hat{\vec{u}}(t)-\hat{\vec{u}}(t)\right)}{(\Oc)^m}^2
      + \hat{h}^2\sum_{i=1}^m\ltwon{\nabla\pdt\left({\Ritz\hat{{u}}_i(t)-\hat{{u}}_i(t)}\right)}{(\Oc)}^2\Bigg\}
\leq
 C(\A,\hat{u})\hat{h}^{2(\ell+1)}.
    \end{split}
  \end{gather}
\end{Lem} 
\begin{Proof} 
Using (\ref{eqn:coercivity}) and (\ref{eqn:ritz_reference_definition}) we have for $i=1,\dotsc,m$
\begin{equation}
\begin{split}
\mu\ltwon{\nabla\hat{\varepsilon_i}}{(\Oc)}^2&
\leq
 a(\hat \varepsilon_i, \hat \Phi -\hat u_i)\Foreach\Phi\in\Vc\\
&
 \leq
 \bar{\mu}\ltwon{\nabla\hat\varepsilon_i}{(\Oc)}\ltwon{\nabla(\Lagrange\hat u_i-\hat u_i)}{(\Oc)}
 \leq
 C\hat h^{\ell}\ltwon{\nabla\hat\varepsilon_i}{(\Oc)}\lv \hat u_i\rv_{\Hil{\ell+1}(\Oc)},
\end{split}
\end{equation}
which shows the energy norm bound of
(\ref{eqn:Ritz_projection_error}). To show the $\Lp{2}$
estimate we use duality. Fix a $t\in(0,T]$ and consider the solution $\hat\psi$ of following elliptic problem
\begin{equation}\label{eqn:dual_problem}
-\nabla\cdot(\vec{B}_t\nabla \hat\psi)=\hat\phi\mbox{ in }\Oc,\quad \vec{B}_t\nabla\hat\psi\cdot\hat{\normal}=0\mbox{ on }\partial\Oc, \quad\int_\Oc\hat{\psi}=0.
\end{equation}
Note that $\ltwon{\hat \psi}{(\Oc)} 
\leq
 C\lv \hat\psi \rv_{\Hil{1}(\Oc)}$ as for any $\hat v$
\begin{equation}
  \inf_{r\in\Reals}\ltwon{\hat v-r}{(\Oc)}= 
  \ltwon[\Big]{\hat v-\frac1{\vert\Oc\vert}\int_{\Oc}\hat v}{(\Oc)} 
\leq
 C\lv\hat v\rv_{\Hil{1}(\Oc)}.
\end{equation}
We therefore have
\begin{equation}\label{eqn:grad_psi_phi}
\mu\ltwon{\nabla\hat \psi}{(\Oc)}^2
\leq
 \bbil{\hat\psi}{\hat\psi}{\Oc}=\ltwop{\hat \phi}{\hat \psi}{\Oc}
\leq
 C\ltwon{\hat\phi}{(\Oc)}\ltwon{\nabla\hat\psi}{(\Oc)}.
\end{equation}
Furthermore we have the estimate
\begin{equation}
\lv\hat\psi\rv_{\Hil{2}(\Oc)}
\leq
 C\ltwon{\lap\hat\psi}{(\Oc)}
\leq
 C\ltwon{\vec{B}\lap\hat\psi}{(\Oc)}=C\ltwon{\hat\phi+\nabla\cdot\vec{B}\cdot\nabla\hat\psi}{(\Oc)}
\leq
 C\ltwon{\hat\phi}{(\Oc)}.
\end{equation}
\changes{Here we have introduced the notation that the divergence of the tensor $\vec{B}$ is a vector defined such that for $i=1,\dotsc,d$,
$(\nabla\cdot\vec B)_i=\sum_{j=1}^d\partial_{x_j}\vec B_{i,j}.$
}
Thus testing (\ref{eqn:dual_problem}) with $\hat\varepsilon_i$ and using (\ref{eqn:ritz_reference_definition}) we have
\begin{equation}\begin{split}
\ltwop{\hat\varepsilon_i}{\hat\phi}{\Oc}&=\bbil{\hat\varepsilon_i}{\hat\psi-\hat\Phi}{\Oc}\Foreach\hat\Phi\in\Vc\\
&
\leq
\bar{	\mu}\ltwon{\nabla\hat\varepsilon_i}{(\Oc)}\ltwon{\nabla(\hat\psi-\Lagrange\hat\psi)}{(\Oc)}
\leq
 C\hat h^{\ell+1}\lv \hat u\rv_{\Hil{\ell+1}(\Oc)}\lv \psi \rv_{\Hil{2}(\Oc)}
\leq
 C\hat h^{\ell+1},
\end{split}
\end{equation}
which completes the proof of
(\ref{eqn:Ritz_projection_error}).  For the proof of
(\ref{eqn:dt_Ritz_projection_error})  using
(\ref{eqn:ritz_dt_defn}) and the fact that the gradient commutes with the time derivative (as we work on the reference domain) we have that \fim and for each $\hat{\Phi}\in\Vc$,
\begin{equation}
\begin{split}\label{eqn:dt_energy_Ritz_I}
\mu&\ltwon{\nabla\pdt\hat\varepsilon_i}{(\Oc)}^2
\leq
 \bbil{\pdt\hat\varepsilon_i}{\pdt\hat\varepsilon_i}{\Oc}= \bbil{\pdt\hat\varepsilon_i}{\hat{\Phi}-\pdt\hat{u}_i}{\Oc}+\bbil{\pdt\hat\varepsilon_i}{\pdt\Ritz\hat{u}_i-\hat{\Phi}}{\Oc}\\
&=\bbil{\pdt\hat\varepsilon_i}{\hat\Phi-\pdt\hat u_i}{\Oc}+\ltwop{\pdt\vec B\nabla\hat\varepsilon_i}{\nabla(\hat\Phi-\pdt\Ritz\hat u_i)}{\Oc}.
\end{split}
\end{equation}
Taking $\hat\Phi=\Lagrange\pdt\hat{u}_i$ in (\ref{eqn:dt_energy_Ritz_I}) gives
\begin{equation}
\begin{split}
\mu\ltwon{\nabla\pdt\hat\varepsilon_i}{(\Oc)}^2
\leq
& C\hat h^\ell\lv \pdt\hat u_i\rv_{\Hil{\ell +1}(\Oc)}\ltwon{\nabla\pdt\hat\varepsilon_i}{(\Oc)}\\
&+\linfn{\pdt\vec B}{(\Oc)}\ltwon{\nabla\hat\varepsilon_i}{(\Oc)}\bigg(\ltwon{\nabla\pdt\hat\varepsilon_i}{(\Oc)}+\big\|{\nabla\left(\Lagrange\pdt\hat u_i-\pdt\hat u_i\right)}\big\|_{\Lp{2}(\Oc)}\bigg)\\
\leq
& \frac{\mu}{2}\|{\nabla\pdt\hat\varepsilon_i}\|_{\Lp{2}(\Oc)}^2+C(\ltwon{\nabla\hat\varepsilon_i}{(\Oc)}^2+\hat h^{2\ell}\lv\pdt \hat u_i\rv_{\Hil{\ell+1}(\Oc)}^2),
\end{split}
\end{equation}
where we have use Young's inequality in the final step.  The previous
estimate (\ref{eqn:Ritz_projection_error}) completes the
proof of the energy norm bound in
(\ref{eqn:dt_Ritz_projection_error}).  For the $\Lp{2}$
estimate we once again use duality. Testing problem
(\ref{eqn:dual_problem}) with $\pdt\hat\varepsilon_i$ and using
(\ref{eqn:ritz_dt_defn}), we have \fim and any
$\hat\Phi\in\Vc$
\begin{equation}\label{eqn:dt_ritz_error_1}
\begin{split}
\ltwop{\pdt\hat\varepsilon_i}{\hat \phi}{\Oc}&=\bbil{\pdt\hat\varepsilon_i}{\hat\psi-\hat\Phi}{\Oc}-\ltwop{(\pdt\vec B)\nabla\hat\varepsilon_i}{\nabla\hat{\Phi}}{\Oc}\\
&=\bbil{\pdt\hat\varepsilon_i}{\hat\psi-\hat\Phi}{\Oc}+\ltwop{(\pdt\vec B)\nabla\hat\varepsilon_i}{\nabla\left(\hat\psi-\hat{\Phi}\right)}{\Oc}-\ltwop{(\pdt\vec B)\nabla\hat\varepsilon_i}{\nabla\hat\psi\rangle}{\Oc}.
\end{split}
\end{equation}
Taking $\hat\Phi=\Lagrange\pdt\hat{u}_i$ in (\ref{eqn:dt_ritz_error_1}) gives
\begin{equation}
\begin{split}
&\lv\ltwop{\pdt\hat\varepsilon_i}{\hat\phi}{\Oc}\rv
\leq
 C\lv\hat\psi\rv_{\Hil{2}(\Oc)}\Bigg(\hat h\bar{\mu}\ltwon{\nabla\pdt\hat\varepsilon_i}{(\Oc)}+\hat h\linfn{\pdt\vec{B}}{(\Oc)}\ltwon{\nabla\hat\varepsilon_i}{(\Oc)}\\
 &\quad+\linfn{\pdt\vec{B}}{\Oc}\ltwon{\varepsilon_i}{(\Oc)}\Bigg)+\ltwon{\nabla\hat\psi}{(\Oc)}\linfn{\nabla\pdt\vec{B}}{(\Oc)}\ltwon{\varepsilon_i}{(\Oc)},
\end{split}
\end{equation} 
where we have used integration by parts to estimate the last term in (\ref{eqn:dt_ritz_error_1}). The previous estimates  and Assumption \ref{assum:mapping} complete the proof.
\end{Proof}
\begin{The}[A priori estimate for the semidiscrete scheme]\label{thm:semidiscrete}
 Suppose Assumptions \ref{assum:reaction} and \ref{assum:well_posed} hold. Furthermore, let Assumption \ref{assum:mapping} hold (with $k=\ell$).  Finally let $\hat{\vec{u}}^h$ be the solution to Problem (\ref{eqn:sd_weak_form}).  Then, the following optimal a priori error estimate holds for the error in the semidiscrete scheme:
\begin{equation}
\begin{split}\label{eqn:semi_discrete_apriori_estimate}
\sup_{t\in[0,T]}\left\{\ltwon{\hat{\vec{u}}^h(t)-\hat{\vec{u}}(t)}{(\Oc)^m}^2\right\}+\sum_{i=1}^m\int_{0}^T\hat{h}^2\ltwon{\nabla(\hat{u}_i^h(t)-\hat{u}_i(t))}{(\Oc)^m}^2\dif t\leq{C}\left(\A,\hat{u}, \widetilde{C}\right)\hat{h}^{2(\ell+1)},
\end{split}
\end{equation}
\end{The}
\begin{Proof} 
Using the decomposition (\ref{eqn:error_decomp_semi_disc}) and Lemma \ref{lem:epsilon_grad_epsilon} we have a bound on the elliptic error and it simply remains to estimate the parabolic error $\hat{\vec{\rho}}^h$. To this end, we use (\ref{eqn:sd_weak_form}) to construct a PDE for $\hat{\rho}^h_i$ by  inserting $\hat{\rho}^h_i$ in place of  $\hat{u}^h_i$ and taking $\hat{\Phi} = \hat{\rho}^h_i$. Using (\ref{eqn:gradient_evolving_eulerian_correspondence}) we obtain \fim 
\begin{equation}\label{eqn:SD_convergence_weak_form_I}
\ltwop{\pdt\left(J\hat{\rho}^h_i\right)}{\hat{\rho}^h_i}{\Oc}+D_i\seminorm{\nabla\hat{\rho}^h_i}{B}^2=\ltwop{\widetilde{f}_i(\hat{\vec{u}}^h)}{J\hat{\rho}^h_i}{\Oc}-\ltwop{\pdt\left(J\Ritz\hat{u}_i\right)}{\hat{\rho}^h_i}{\Oc}-\bbil{\Ritz\hat{u}_i}{\hat{\rho}^h_i}{\Oc}.
\end{equation}
Using  (\ref{eqn:cont_weak_form_reference_ddtju}), (\ref{eqn:f_tilde_equals_f_for_exact_solution}) and (\ref{eqn:ritz_reference_definition}) gives
\begin{equation}\label{eqn:SD_convergence_weak_form_II}
\ltwop{\pdt\left(J\hat{\rho}^h_i\right)}{\hat{\rho}^h_i}{\Oc}+D_i\seminorm{\nabla\hat{\rho}^h_i}{B}^2=\ltwop{\widetilde{f}_i(\hat{\vec{u}}^h)-\widetilde{f}_i(\hat{\vec{u}})}{J\hat{\rho}^h_i}{\Oc}-\ltwop{\pdt\left(J\hat{\varepsilon}_i\right)}{\hat{\rho}^h_i}{\Oc}.
\end{equation}
Dealing with the first term on the left of (\ref{eqn:SD_convergence_weak_form_II}) as in (\ref{eqn:pdtJhatu_u_on_evolving_frame}):
\begin{equation}\label{eqn:pdtJhatrho_rho_on_evolving_frame}
\ltwop{\pdt\left(J\hat{\rho}^h_i\right)}{\hat{\rho}^h_i}{\Oc}=\frac{1}{2}\left(\frac{\dif}{\dif t}\mynorm{\hat{\rho}_i^h}{J}^2+\ltwop{J\hat{\rho}^h_i\nabla\cdot\vec{a}(\A_t(\vec \xi)))}{\hat{\rho}^h_i}{\Oc}\right).
\end{equation}
Dealing with the first term on the right of (\ref{eqn:SD_convergence_weak_form_II}) using (\ref{eqn:error_decomp_semi_disc}) and the MVT we have
\begin{equation}
\lv\ltwop{\widetilde{f}_i(\hat{\vec{u}}^h)-\widetilde{f}_i(\hat{\vec{u}})}{J\hat{\rho}^h_i}{\Oc}\rv\leq\widetilde{C}\left(\ltwop{\sum_{j=1}^m\left(\lv\hat{\varepsilon}_j\rv+\lv\hat{\rho}^h_j\rv\right)}{J\lv\hat{\rho}^h_i\rv}{\Oc}\right).
\end{equation}
Applying Young's inequality:
\begin{equation}
\begin{split}
\lv\ltwop{\widetilde{f}_i(\hat{\vec{u}}^h)-\widetilde{f}_i(\hat{\vec{u}})}{J\hat{\rho}^h_i}{\Oc}\rv\leq\widetilde{C}\Bigg(\left(m+\frac{1}{2}\right)\mynorm{\hat{\rho}^h_i}{J}^2+\sum_{j\not=i}\frac{1}{2}\mynorm{\hat{\rho}^h_j}{J}^2+\frac{1}{2}\mynorm{\vec{\hat{\varepsilon}}}{J^m}^2\Bigg).
\end{split}
\end{equation}
Summing over $i$ we have
\begin{equation}\label{eqn:sd_estimate_R1}
\begin{split}
\sum_{i=1}^m\lv\ltwop{\widetilde{f}_i(\hat{\vec{u}}^h)-\widetilde{f}_i(\hat{\vec{u}})}{J\hat{\rho}^h_i}{\Oc}\rv\leq\widetilde{C}\Bigg(\frac{3m}{2}\mynorm{\hat{\vec{\rho}}^h}{J^m}^2+\frac{m}{2}\mynorm{\vec{\hat{\varepsilon}}}{J^m}^2\Bigg).
\end{split}
\end{equation}
Dealing with the second term on the right of (\ref{eqn:SD_convergence_weak_form_II}):
\begin{equation}
\begin{split}
\lv\ltwop{\pdt\left(J\hat{\varepsilon}_i\right)}{\hat{\rho}^h_i}{\Oc}\rv\leq&\lv\ltwop{J\pdt\hat{\varepsilon}_i}{\hat{\rho}^h_i}{\Oc}\rv+\lv\ltwop{\pdt\left(J\right)\hat{\varepsilon}_i}{\hat{\rho}^h_i}{\Oc}\rv\\
\leq&\frac{1}{2}\Bigg(\mynorm{\hat{\rho}^h}{J}^2+\ltwop{{J}\pdt\hat{\varepsilon}_i}{\pdt\hat{\varepsilon}_i}{\Oc}+\ltwop{\lv\pdt({J})\rv\hat{\rho}^h_i}{\hat{\rho}^h_i}{\Oc}+\ltwop{\lv\pdt({J})\rv\hat{\varepsilon}_i}{\hat{\varepsilon}_i}{\Oc}\Bigg),
\end{split}
\end{equation}
where we have used Young's inequality for the second step. Now using (\ref{eqn:dtdetj}) and summing over $i$ we have
\begin{equation}
\begin{split}\label{eqn:sd_estimate_R2}
\sum_{i=1}^m\lv\ltwop{\pdt\left(J\hat{\varepsilon}_i\right)}{\hat{\rho}^h_i}{\Oc}\rv\leq\frac{1}{2}\Bigg(&\mynorm{\hat{\vec{\rho}}^h}{J^m}^2+\ltwop{J\hat{\vec{\rho}^h}\lv\nabla\cdot\vec{a}\left(\A_t(\vec{\xi}),t\right)\rv}{\hat{\vec{\rho}}^h}{\Oc^m}\\
&+\mynorm{\pdt\vec{\hat{\varepsilon}}}{J^m}^2+\linfn{\pdt{J}}{(\Oc\times[0,T])}\ltwon{\vec{\hat{\varepsilon}}}{(\Oc)^m}^2\Bigg).
\end{split}
\end{equation}
Combining (\ref{eqn:pdtJhatrho_rho_on_evolving_frame}), (\ref{eqn:sd_estimate_R1}), (\ref{eqn:sd_estimate_R2}) 
\begin{equation}
\begin{split}
\frac{\dif}{\dif t}\mynorm{\hat{\vec{\rho}}^h}{J^m}^2+2\sum_{i=1}^m D_i\seminorm{\nabla\hat{\rho}^h_i}{B}^2\leq{C}\Bigg(\mynorm{\hat{\vec{\rho}}^h}{J^m}^2+\ltwon{\vec{\hat{\varepsilon}}}{(\Oc)^m}^2+\ltwon{\pdt\vec{\hat{\varepsilon}}}{(\Oc)^m}^2\Bigg),
\end{split}
\end{equation}
where we have used the fact that Assumption \ref{assum:mapping} implies $J,\pdt{J}\in\Lp{\infty}{\left(\smash{\Oc}\times[0,T]\right)}$. 
Integrating in time, using Lemma \ref{lem:epsilon_grad_epsilon} and applying Gronwall's Lemma we have
\begin{equation}
\mynorm{\hat{\vec{\rho}}^h(t)}{J^m}^2+2\sum_{i=1}^m D_i\int_{0}^T\seminorm{\nabla\hat{\rho}^h_i}{B}^2\leq{C}\left(\mynorm{\hat{\vec{\rho}}^h(0)}{J^m}^2+\hat{h}^{2(\ell+1)}\right).
\end{equation}
To estimate $\hat{\vec{\rho}}^h(0)$, we note
\begin{equation}
\begin{split}\label{eqn:rho_zero_semidiscrete_estimate}
\mynorm{\hat{\vec{\rho}}^h(0)}{J^m}^2\leq\mynorm{\hat{\vec{u}}(0)-\Lagrange\hat{\vec{u}}(0)}{J^m}+\mynorm{\hat{\vec{\varepsilon}}^h}{J^m}
\leq{C}\hat{h}^{\ell+1},
\end{split}
\end{equation}
where we have used (\ref{eqn:reference_fe_space_accuracy}), the assumption on the regularity of the exact solution and Lemma \ref{lem:epsilon_grad_epsilon} in the last step. Assumption \ref{assum:mapping} and the equivalence of norms  (\ref{eqn:norm_evolving_eulerian_correspondance}) completes the proof.
\end{Proof}
\section{Error analysis of the fully discrete approximation}\label{subsecn.:fully_disc}
In this section we provide the convergence result for the fully
discrete scheme (\ref{eqn:fully_discrete_scheme}).  The main result of
this paper is Theorem \ref{thm:fully_discrete}, whose proof is given in detail below. 
We follow that up with a convergence result in the
$\Lp\infty(\Oc)$ norm which allows the use of the original $\vec f$
(without extending to $\widetilde{\vec f}$ in the numerical method).
\begin{The}[A priori estimate for the fully discrete scheme]
  \label{thm:fully_discrete}
  Suppose Assumptions \ref{assum:reaction} and \ref{assum:well_posed}
  hold. Suppose Assumption \ref{assum:mapping} (with
  $k=\ell$) holds.  Let $\hat{\vec{U}}$ be the solution to
  (\ref{eqn:fully_discrete_scheme}).  Suppose the timestep
  satisfies a stability condition defined in
  (\ref{eqn:timestep_restriction}). Then, the following optimal a
  priori estimate holds for the error in the fully discrete scheme:
  \begin{equation}
  \begin{split}
    \label{eqn:fully_discrete_apriori_estimate}
    \ltwon{\hat{\vec{U}}^n-\hat{\vec{u}}^n}{(\Oc)^m}^2+&\tau\hat{h}^2\sum_{i=1}^m{D}_i\ltwon{\nabla\left(\hat{{U}}_i^n-\hat{{u}}_i^n\right)}{(\Oc)}^2\\
    &\leq{C}\left(\A,\hat{u}, \widetilde{C}\right)\left(\hat{h}^{2(\ell+1)}+\tau^2\right), \quad \mbox{for } n\in[0,\dotsc,N],
\end{split}
 \end{equation}
      with $\widetilde{C}$ as defined in (\ref{eqn:ftilde_def}).
\end{The}
\additions{
\begin{Rem}[Error estimate for the evolving domain scheme]
The schemes (\ref{eqn:fully_discrete_scheme}) and (\ref{eqn:moving_fully_discrete_scheme}) are equivalent. Thus Theorem \ref{thm:fully_discrete} also provides an error estimate for the evolving domain based scheme (\ref{eqn:moving_fully_discrete_scheme}).
\end{Rem}
}
\begin{Proof}[of Theorem \ref{thm:fully_discrete}] Decomposing the error as in (\ref{eqn:error_decomp_semi_disc}) we have
\begin{equation}\label{eqn:error_decomposition_fully_discrete}
\begin{split}
\ltwon{\hat{\vec{U}}^n-\hat{\vec{u}}^n}{(\Oc)^m}^2&\leq\ltwon{\Ritz\hat{\vec{u}}^n-\hat{\vec{u}}^n}{(\Oc)^m}^2+\ltwon{\hat{\vec{U}}^n-\Ritz\hat{\vec{u}}^n}{(\Oc)^m}^2=\ltwon{\hat{\vec{\varepsilon}}^n}{(\Oc)^m}^2+\ltwon{\hat{\vec{\rho}}^n}{(\Oc)^m}^2.
\end{split}
\end{equation}
From Lemma \ref{lem:epsilon_grad_epsilon} we have the following bound on the elliptic error: 
\begin{equation}
\ltwon{\hat{\vec{\varepsilon}}^n}{(\Oc)^m}^2\leq{C}\hat{h}^{2(\ell+1)} \quad \text{ for }n\in[0,\dotsc,N].
\end{equation}
Therefore it only remains to estimate $\hat{\vec{\rho}}^n$. Constructing an expression for $\hat{\vec{\rho}}^n$ as in (\ref{eqn:SD_convergence_weak_form_I}), using (\ref{eqn:fully_discrete_scheme}) and (\ref{eqn:ritz_reference_definition}) we obtain \fim
\begin{equation}\label{eqn:fd_rho_weak_form}
\begin{split}
&\ltwop{\bar{\partial}[J\hat{\rho}_i]^n}{\hat{\rho}_i^n}{\Oc}+D_i\seminorm{\nabla\hat{\rho}^n_i}{B}^2=\ltwop{\hat{U}_i^n\widetilde{F}_i(\hat{\vec{U}}^{n-1})}{[J\hat{\rho}_i]^n}{\Oc}\\
&\qquad-\ltwop{\bar{\partial}[J\Ritz\hat{u}_i]^n}{\hat{\rho}_i^n}{\Oc}-D_i\ltwop{[J\vec{K}\nabla\hat{u}_i]^n}{[\vec{K}\nabla\hat{\rho}_i]^n}{\Oc}\\
&\quad=\ltwop{\hat{U}_i^n\widetilde{F}_i(\hat{\vec{U}}^{n-1})-\widetilde{f}_i(\hat{\vec{u}}^n)}{[J\hat{\rho}_i]^n}{\Oc}-\ltwop{\bar{\partial}[J\hat{\varepsilon}_i]^n}{\hat{\rho}_i^n}{\Oc}+\ltwop{\left(\bar{\partial}-\pdt\right)[J\hat{u}_i]^n}{\hat{\rho}_i^n}{\Oc},
\end{split}
\end{equation}
where we have used (\ref{eqn:cont_weak_form_reference_ddtju}) for the second step and $\vec{\widetilde{F}}$ is as defined in (\ref{eqn:ftilde_def}).
Using Young's inequality for the first term on the left hand side of (\ref{eqn:fd_rho_weak_form}) gives 
 \begin{equation}
 \begin{split}
 \ltwop{\bar{\partial}[J\hat{\rho}_i]^n}{\hat{\rho}_i^n}{\Oc}
 \geq\frac{1}{\tau_n}\Bigg(&\mynorm{\hat{\rho}^n_i}{J}^2
 -
 \frac{1}{2}\left(\ltwop{J^{n-1}\hat{\rho}_i^n}{\hat{\rho}_i^n}{\Oc}
 +
 \ltwop{J^{n-1}\hat{\rho}_i^{n-1}}{\hat{\rho}_i^{n-1}}{\Oc}\right)\Bigg),
\end{split}
\end{equation}
where we have used (\ref{eqn:norm_evolving_eulerian_correspondance}). Summing over $i$ we have
 \begin{equation}\label{eqn:dbarJhatrho_rho_on_evolving_frame}
 \begin{split}
 \sum_{i=1}^m\ltwop{\bar{\partial}[J\hat{\rho}_i]^n}{\hat{\rho}_i^n}{\Oc}\geq\frac{1}{\tau_n}\left(1-\frac{1}{2}\linfn{\frac{J^{n-1}}{J^{n}}}{(\Oc)}\right)\mynorm{\hat{\vec{\rho}}^n}{J^m}^2-\frac{1}{2\tau_n}\mynorm{\hat{\vec{{\rho}}}^{n-1}}{J^m}^2.
\end{split}
\end{equation}
Using \ref{eqn:error_decomposition_fully_discrete} and the MVT for the first term on the right hand side of (\ref{eqn:fd_rho_weak_form}) gives 
\begin{equation}
\begin{split}
\Bigg\vert&\ltwop{\hat{U}_i^n\widetilde{F}_i(\hat{\vec{U}}^{n-1})-\widetilde{f}_i(\hat{\vec{u}}^n)}{[J\hat{\rho}_i]^n}{\Oc}\Bigg\vert\\
&\quad\leq\widetilde{C}\sum_{j=1}^m\ltwop{\lv\hat{\varepsilon}^{n-1}_j\rv+\lv\hat{\rho}^{n-1}_j\rv+\lv\tau_n\bar{\partial}\hat{u}^{n}_j\rv+\lv\hat{\varepsilon}^{n}_i\rv+\lv\hat{\rho}^{n}_i\rv}{J^n\lv\hat{\rho}_i^n\rv}{\Oc}\\
&\quad\leq{C}\widetilde{C}\Bigg(\mynorm{\hat{\rho}_i^n}{J}^2+\linfn{\frac{J^n}{J^{n-1}}}{(\Oc)}\mynorm{\hat{\vec{\rho}}^{n-1}}{J^m}^2\\
&\qquad+\linfn{J^n}{(\Oc)}\bigg(\ltwon{\hat{{\varepsilon}}_i^n}{(\Oc)}^2+\ltwon{\hat{\vec{\varepsilon}}^{n-1}}{(\Oc)^m}^2+\ltwon{\tau_n\bar{\partial}\hat{\vec{u}}^n}{(\Oc)^m}^2\bigg)\Bigg).
\end{split}
\end{equation}
where we have used Young's inequality for the second step. Summing over $i$ we have
\begin{equation}\label{eqn:fd_estimate_R1}
\begin{split}
\sum_{i=1}^m&\lv\ltwop{\hat{U}_i^n\widetilde{F}_i(\hat{\vec{U}}^{n-1})-\widetilde{f}_i(\hat{\vec{u}}^n)}{[J\hat{\rho}_i]^n}{\Oc}\rv\leq{C}\widetilde{C}\Bigg(\mynorm{\hat{\vec{\rho}}^n}{J^m}^2+\linfn{\frac{J^n}{J^{n-1}}}{(\Oc)}\mynorm{\hat{\vec{\rho}}^{n-1}}{J^m}^2\\
&+\linfn{J^n}{(\Oc)}\bigg(\ltwon{\hat{\vec{\varepsilon}}^n}{(\Oc)^m}^2+\ltwon{\hat{\vec{\varepsilon}}^{n-1}}{(\Oc)^m}^2+\ltwon{\tau_n\bar{\partial}\hat{\vec{u}}^n}{(\Oc)^m}^2\bigg)\Bigg),
\end{split}
\end{equation}
Applying Young's inequality to the second and third term on the right of (\ref{eqn:fd_rho_weak_form}) gives
\begin{equation}\label{eqn:fd_estimate_R23}
\begin{split}
\lv\ltwop{\bar{\partial}[J\hat{\varepsilon}_i]^n}{\hat{\rho}_i^n}{\Oc}\rv&+\lv\ltwop{\left(\bar{\partial}-\pdt\right)[J\hat{u}_i]^n}{\hat{\rho}_i^n}{\Oc}\rv\\
\leq\mynorm{\hat{\rho}^n_i}{J}^2&+\frac{1}{2}\linfn{\frac{1}{J^{n}}}{(\Oc)}\left(\ltwon{\bar{\partial}[J\hat{\varepsilon}_i]^n}{(\Oc)}^2+\ltwon{\left(\bar{\partial}-\pdt\right)[J\hat{u}_i]^n}{(\Oc)}^2\right).
\end{split}
\end{equation}
Using (\ref{eqn:dbarJhatrho_rho_on_evolving_frame}), (\ref{eqn:fd_estimate_R1}) and (\ref{eqn:fd_estimate_R23}) in (\ref{eqn:fd_rho_weak_form}) gives
\begin{equation}
\begin{split}
\frac{1}{\tau_n}&\left(1-\frac{1}{2}\linfn{\frac{J^{n-1}}{J^{n}}}{(\Oc)}-C\widetilde{C}\tau_n\right)\ltwon{\hat{\vec{\rho}}^n}{J^m}^2+\sum_{i=1}^m D_i\seminorm{\nabla\hat{\rho}^n_i}{B}^2\\
&\leq\left(\frac{1}{2\tau_n}+C\widetilde{C}\linfn{\frac{J^n}{J^{n-1}}}{(\Oc)}\right)\mynorm{\hat{\vec{{\rho}}}^{n-1}}{J^m}^2+{C}\widetilde{C}\linfn{J^n}{(\Oc)}\Big(\ltwon{\hat{\vec{\varepsilon}}^n}{(\Oc)^m}^2+\ltwon{\hat{\vec{\varepsilon}}^{n-1}}{(\Oc)^m}^2\\
&+\ltwon{\tau_n\bar{\partial}\hat{\vec{u}}^n}{(\Oc)^m}^2\Big)+\frac{1}{2}\linfn{\frac{1}{J^{n}}}{(\Oc)}\Big(\ltwon{\bar{\partial}[J\hat{\vec{\varepsilon}}]^n}{(\Oc)^m}^2+\ltwon{\left(\bar{\partial}-\pdt\right)[J\hat{\vec{u}}]^n}{(\Oc)^m}^2\Big).
\end{split}
\end{equation}
Let $\tau^\prime>0$ be such that, for $\tau<\tau^\prime$ and for $n=1,\dotsc,N$,
\begin{equation}\label{eqn:timestep_restriction}
1-\frac{1}{2}\linfn{\frac{J^{n-1}}{J^{n}}}{(\Oc)}-C\widetilde{C}\tau>0.
\end{equation}
Such a $\tau^{\prime}$ exists since
\begin{equation}
\lim_{\tau\to0}\left\{\frac{1}{2}\linfn{\frac{J^{n-1}}{J^{n}}}{(\Oc)}+C\widetilde{C}\tau\right\}=\frac{1}{2}.
\end{equation}
For $\tau<\tau^\prime$, we have
\begin{equation}\label{eqn:rho_n_expression_rho_zero_residual_n}
\mynorm{\hat{\vec{\rho}}^n}{J^m}^2+\sum_{i=1}^m C\tau{D}_i\seminorm{\nabla\hat{\rho}^n_i}{B}^2\leq{C}\left(\bar{C}^n\mynorm{\hat{\vec{{\rho}}}^{n-1}}{J^m}^2+\tau\mathcal{R}^n\right),
\end{equation}
where $\bar{C}^n=1+\tau\widetilde{C}\linfn{\frac{J^n}{J^{n-1}}}{(\Oc)}$ and 
\begin{equation}\label{eqn:residual_definition}
\begin{split}
\mathcal{R}^n:=\widetilde{C}&\linfn{J^n}{(\Oc)}\Big(\ltwon{\hat{\vec{\varepsilon}}^n}{(\Oc)^m}^2+\ltwon{\hat{\vec{\varepsilon}}^{n-1}}{(\Oc)^m}^2+\ltwon{\tau\bar{\partial}\hat{\vec{u}}^n}{(\Oc)^m}^2\Big)\\
+\frac{1}{2}&\linfn{\frac{1}{J^{n}}}{(\Oc)}\Big(\ltwon{\bar{\partial}[J\hat{\vec{\varepsilon}}]^n}{(\Oc)^m}^2+\ltwon{\left(\bar{\partial}-\pdt\right)[J\hat{\vec{u}}]^n}{(\Oc)^m}^2\Big).
\end{split}
\end{equation}
Therefore, for $n=1,\dotsc,N$,
\begin{equation}
\begin{split}
\mynorm{\hat{\vec{\rho}}^n}{J^m}^2+\sum_{i=1}^m C\tau{D}_i\seminorm{\nabla\hat{\rho}^n_i}{B}^2\leq{C}\Biggr(\prod_{k=1}^n\bar{C}^k\mynorm{\hat{\vec{{\rho}}}^{0}}{J^m}^2+\tau\sum_{j=1}^n\prod_{i=j}^n\bar{C}^i\mathcal{R}^j\Biggr).
\end{split}
\end{equation}
\additions{
 For $n=1,\dotsc,N$, we have
\begin{equation}
\overline{C}^n=1+\tau\widetilde{C}\left\|\frac{J^n}{J^{n-1}}\right\|_{L^{\infty}(\hat\Omega)}
\leq 1+\tau\widetilde{C}\left\|{J^n}\right\|_{L^{\infty}(\hat\Omega)}\left\|\frac{1}{J^{n-1}}\right\|_{L^{\infty}(\hat\Omega)}
\leq 1+\tau \widetilde{C}C,
\end{equation}
where the last line follows by Assumption \ref{assum:mapping}. Thus  $0< \Pi_{i=j}^n \overline{C}^i\leq\Pi_{k=1}^n \overline{C}^k\leq\left(1+\tau \widetilde{C}C\right)^n$.}

Considering the first two terms on the right of (\ref{eqn:residual_definition}), we have for $n=1,\dotsc,N$
\begin{equation}
\begin{split}\label{eqn:residual_1_2}
\widetilde{C}\linfn{J^n}{(\Oc)}\Big(\ltwon{\hat{\vec{\varepsilon}}^n}{(\Oc)^m}^2+\ltwon{\hat{\vec{\varepsilon}}^{n-1}}{(\Oc)^m}^2\Big)\leq2\widetilde{C}\sup_{s\in[0,\dotsc,N]}\linfn{J^s}{(\Oc)}\ltwon{\hat{\vec{\varepsilon}}^{s}}{(\Oc)^m}^2\leq\widetilde{C}{C}\hat{h}^{2(\ell+1)},
\end{split}
\end{equation}
where we have used Assumption \ref{assum:mapping} and Lemma \ref{lem:epsilon_grad_epsilon}.
Dealing with the third term on the right of (\ref{eqn:residual_definition}), we have
\begin{equation}
\begin{split}\label{eqn:residual_3}
\widetilde{C}\linfn{J^n}{(\Oc)}\ltwon{\tau\bar{\partial}\hat{\vec{u}}^n}{(\Oc)^m}^2=\widetilde{C}\linfn{J^n}{(\Oc)}\ltwon{\int_{t^{n-1}}^{t^{n}}\pdt\hat{\vec{u}}^s\dif s}{(\Oc)^m}^2\leq\widetilde{C}{C}\tau^{2},
\end{split}
\end{equation}
where we have used Assumptions \ref{assum:well_posed} and \ref{assum:mapping}.
For the fourth term on the right of (\ref{eqn:residual_definition}) we have
\begin{equation}
\begin{split}\label{eqn:residual_4}
\frac{1}{2}\linfn{\frac{1}{J^{n}}}{(\Oc)}\ltwon{\bar{\partial}[J\hat{\vec{\varepsilon}}]^n}{(\Oc)^m}^2&\leq
\frac{1}{2}\linfn{\frac{1}{J^{n}}}{(\Oc)}\ltwon{\frac{1}{\tau_n}{\int_{t^{n-1}}^{t^{n}}}\pdt[J\hat{\vec{\varepsilon}}]^s\dif s}{(\Oc)^m}^2\\
&\leq{C}\sup_{s\in[t^{n-1},t^{n}]}\ltwon{\hat{\vec{\varepsilon}}^s}{(\Oc)^m}^2\leq{C}\hat{h}^{2(\ell+1)},
\end{split}
\end{equation}
where we have used Assumption \ref{assum:mapping} for the second step and Lemma \ref{lem:epsilon_grad_epsilon} for the final step. Finally, for the fifth term on the right of (\ref{eqn:residual_definition}) we have 
\begin{equation}
\begin{split}\label{eqn:residual_5}
\linfn{\frac{1}{J^{n}}}{(\Oc)}&\ltwon{\left(\bar{\partial}-\pdt\right)[J\hat{\vec{u}}]^n}{(\Oc)^m}^2=\linfn{\frac{1}{J^{n}}}{(\Oc)}\ltwon{\frac{1}{\tau_n}\int^{t^n}_{t^{n-1}}\left(s-t^{n-1}\right)\partial_{t t}[J\hat{\vec{u}}]^s\dif s}{(\Oc)^m}^2\\
&\leq{C}\tau^2\sup_{s\in[t^{n-1},t^n]}\left(\ltwon{\pdt\hat{\vec{u}}^s}{(\Oc)^m}^2+\ltwon{\hat{\vec{u}}^s}{(\Oc)^m}^2\right),
\end{split}
\end{equation}
where we have used Assumption \ref{assum:mapping} for the second step and Assumption \ref{assum:well_posed} for the final step. Combining (\ref{eqn:residual_1_2}), (\ref{eqn:residual_3}), (\ref{eqn:residual_4}) and (\ref{eqn:residual_5}) we have
\begin{equation}\label{eqn:residual_full_estimate}
\mathcal{R}^n\leq{C}\left(\hat{h}^{2(\ell+1)}+\tau^2\right) \text{ for } n=1,\dotsc,N.
\end{equation}
Using  (\ref{eqn:rho_zero_semidiscrete_estimate}) we have
\begin{equation}
\label{eqn:rho_zero_estimate_fully_discrete}
\mynorm{\hat{\vec{\rho}}^0}{J^m}^2=\mynorm{\hat{\vec{\rho}}^h(0)}{J^m}^2\leq{C}\hat{h}^{2(\ell+1)}.
\end{equation}
Applying estimates (\ref{eqn:residual_full_estimate}) and (\ref{eqn:rho_zero_estimate_fully_discrete})  in (\ref{eqn:rho_n_expression_rho_zero_residual_n}) completes the proof of Theorem \ref{thm:fully_discrete}.
\end{Proof}
\begin{Rem}[Stability of the fully discrete scheme]
The timestep restriction (\ref{eqn:timestep_restriction}) is composed
of a term arising from domain growth (the term involving the
determinant $J$ of the diffeomorphism $\A$) and a term arising from
the nonlinear reaction kinetics (the term containing
$\widetilde{C}$). It is worth noting that for a given set of reaction
kinetics, i.e., a given $\widetilde{C}$, larger timesteps are
admissible on growing domains (as we have
$\linfn{\fracl{J^{n-1}}{J^{n}}}{(\Oc)}<1$ for all $n=1,\dotsc,N$). If we consider for
illustrative purposes the heat equation, i.e, the case
$\widetilde{C}=0$, we recover unconditional stability on growing
domains whereas for contracting domains
(\ref{eqn:timestep_restriction}) implies a stability condition on the
timestep dependent on the growth rate.
\end{Rem}

\deletions{
\begin{Rem}[Qualitative estimates on the exact solution]
In practise only qualitative a priori estimates are generally available for the exact solution and the region ${I}_\delta$ defined in Assumption \ref{assum:well_posed} is  not explicitly known. Thus, we cannot construct the function $\widetilde{\vec{f}}$ defined in (\ref{eqn:ftilde_def}). To this end, we introduce the following assumption to circumvent the construction of $\widetilde{\vec{f}}$.
\end{Rem}
\begin{Assum}[Dimension dependent polynomial degree]\label{assum:dim_dependent_accuracy}
We wish to  invoke estimate (\ref{eqn:reference_fe_space_linf_accuracy}) with a positive power of $\hat{h}$ and thus we require the degree of the finite element space to satisfy
$\ell>\frac{d}{2}-1,$
where $d$ is the spatial dimension. Thus, we require  piecewise linear or higher basis functions for $d\leq{2}$ and at least piecewise quadratics for $d=3$.
\end{Assum}
}

\additions
{
In practice only qualitative a priori estimates are generally available for the exact solution and the region ${I}$ defined in Assumption \ref{assum:reaction} is  not explicitly known.  To this end, we show a maximum-norm bound on the discrete solution to circumvent the construction of $\widetilde{\vec{f}}$.

We wish to  invoke estimate (\ref{eqn:reference_fe_space_linf_accuracy}) with a positive power of $\hat{h}$ and thus we require the degree of the finite element space to satisfy
$\ell>\frac{d}{2}-1,$
where $d$ is the spatial dimension. For any physically relevant domain $(d <4)$  piecewise linear or higher basis functions suffice.
\begin{Rem}[Maximum-norm bound of the discrete solution]
Let the assumptions in Theorem
\ref{thm:fully_discrete} be valid and let the degree of the finite element space satisfy
$\ell>\frac{d}{2}-1,$ where $d$ is the spatial dimension. Then 
\begin{equation}
  \label{eqn:maximum-norm-convergence}
  \linfn{\hat{\vec{u}}^n-\hat{\vec{U}}^n}{(\Oc)^m}\leq{C}\hat{h}^{\ell+1-\frac{d}{2}},  
\end{equation}
and for sufficiently small
mesh-size $\hat{h}$ the discrete solution $\hat{\vec{U}}^n$ to Problem
(\ref{eqn:fully_discrete_scheme}) is in the region ${I}$, defined in Assumption \ref{assum:reaction}, for
all $n\in[0,\dotsc,N]$. Thus, we may replace $\widetilde{\vec{F}}$ in
(\ref{eqn:fully_discrete_scheme}) by $\vec{F}$.
\par
Indeed, for $n\in[0,\dotsc,N]$ we have for $\Clement$ the Cl\'ement interpolant
\begin{align}
\linfn{\hat{\vec{u}}^n-\hat{\vec{U}}^n}{(\Oc)^m}
&
\leq
\linfn{\Clement\hat{\vec{u}}^n-\hat{\vec{U}}^n}{(\Oc)^m}+\linfn{\hat{\vec{u}}^n
  -\Clement\hat{\vec{u}}^n}{(\Oc)^m}.
\end{align}
Using (\ref{eqn:reference_fe_space_linf_accuracy}) 
and (\ref{eqn:fe_space_reference_inv_est})
gives
\begin{equation}
\begin{split}
\linfn{\hat{\vec{u}}^n-\hat{\vec{U}}^n}{(\Oc)^m}
&
\leq{C}\Bigg(
\hat{h}^{-d/2}\bigg(\ltwon{\Clement\hat{\vec{u}}^n-\hat{\vec{u}}^n}{(\Oc)^m}+\ltwon{\hat{\vec{u}}^n-\hat{\vec{U}}^n}{(\Oc)^m}\bigg)
\\
&
\phantom{\leq C\Bigg(}
+
\hat{h}^{\ell+1-d/2}\seminorm{\hat{\vec{u}}^n}{\Hil{\ell+1}{(\Oc)^m}}
\Bigg).
\end{split}
\end{equation}
Error bound (\ref{eqn:maximum-norm-convergence}) now follows from
(\ref{eqn:reference_fe_space_accuracy}) and Theorem
\ref{thm:fully_discrete}.
Thus, if $\hat{h}$ is taken sufficiently small we have
\begin{equation}\label{eqn:linf_closeness}
\sup_{n\integerbetween0N}\linfn{\hat{\vec{u}}^n-\hat{\vec{U}}^n}{(\Oc)^m}\leq\delta,
\end{equation}
for any $\delta\in\mathbb{R}^+$. Therefore, $\hat{\vec{U}}^n\in{I}$ for all $n\in[0,\dotsc,N]$ and thus, $\widetilde{\vec{f}}(\hat{\vec{U}})=\vec{f}(\hat{\vec{U}})$. The following corollary follows immediately.
\end{Rem}
\begin{Cor}[Convergence of a practical finite element method]
  \label{cor:convergence-of-a-practical-FEM}
  Let the assumptions in Theorem \ref{thm:fully_discrete} be valid and
  let the degree of the finite element space satisfy
  $\ell>\frac{d}{2}-1$, where $d$ is the spatial dimension. Then, for
  a sufficiently small mesh-size $\hat{h}$ the scheme
  (\ref{eqn:fully_discrete_scheme}), with $\tilde{\vec F}$ replaced by
  $\vec F$ possesses a unique solution
  $\seqsufromto{\smash{\hat{\vec{U}}}}n0N$.  It satisfies the
  following optimal-rate a priori error estimate:
  \begin{equation}
  \begin{split}
    \label{eqn:practical_fully_discrete_apriori_estimate}
    \ltwon{\hat{\vec{U}}^n-\hat{\vec{u}}^n}{(\Oc)^m}^2+&\tau\hat{h}^2\sum_{i=1}^m{D}_i\ltwon{\nabla\left(\hat{{U}}_i^n-\hat{{u}}_i^n\right)}{(\Oc)}^2\\
    &\leq{C}\left(\A,\hat{u}, \widetilde{C}\right)\left(\hat{h}^{2(\ell+1)}+\tau^2\right), \quad \mbox{for } n\in[0,\dotsc,N],
\end{split}
 \end{equation}
     with $\widetilde{C}$ as defined in (\ref{eqn:ftilde_def}).
\end{Cor}
\begin{Rem}[How small must the mesh-size be?]
  Knowing the the meshsize is ``sufficiently small'' in Corollary~\ref{cor:convergence-of-a-practical-FEM}
  is possible, by verifying that the computed solution remains in the region
  $I$ defined in Assumption \ref{assum:reaction},
\end{Rem}
}
\section{Implementation}\label{secn.:implementation}

In this section we illustrate the implementation of the finite element
scheme with explicit nonlinear reaction functions.  We consider the following widely
studied set of reaction kinetics.
\begin{Defn}[Schnakenberg's ``activator-depleted substrate'' model
\citep{schnakenberg1979simple,gierer72,prigo68}]
We consider the following {\it activator depleted substrate} model, also known as
the Brusselator model in nondimensional form:
\begin{equation}\label{eqn:schnak}
\begin{split}
f_1\left(u_{1}, u_{2}\right)=\gamma\left(a-u_{1}+u_{1}^2u_{2}\right)\text{ and }
f_2\left(u_{1}, u_{2}\right)=\gamma\left(b-u_{1}^2u_{2}\right),
\end{split}
\end{equation}
where $0 < a, b, \gamma < \infty $.
\end{Defn}
\additions{
\begin{Rem}[Applicability of Assumption \ref{assum:reaction}]
 The Schnakenberg reaction kinetics satisfy the structural assumptions on the nonlinear reaction vector field as
\begin{equation}
\begin{split}
f_1\left(u_{1}, u_{2}\right)=\gamma \left(a+u_1F_1\left(u_{1}, u_{2}\right)\right)\text{ and }
f_2\left(u_{1}, u_{2}\right)=\gamma \left(b+u_2F_2\left(u_{1}, u_{2}\right)\right),
\end{split}
\end{equation}
where 
 \begin{equation}
\begin{split}
F_1\left(u_{1}, u_{2}\right)=u_1u_2-1\text{ and }
F_2\left(u_{1}, u_{2}\right)=-u_1^2.
\end{split}
\end{equation}
Clearly $\vec f, \vec F \in C^1(\mathbb{R}^2)$ thus Assumption \ref{assum:reaction} holds for the Schnakenberg kinetics.
\end{Rem}
}

In matrix vector form scheme (\ref{eqn:fully_discrete_scheme}) equipped with kinetics (\ref{eqn:schnak}) and appropriate initial approximations  ${\vec{W}}^0_1, {\vec{W}}^0_2$  is: To solve for ${\vec{W}}^n_1, {\vec{W}}^n_2$, $n=[1,\dotsc,N]$, the linear systems given by
\begin{equation}\label{eqn_matrix_vector_scheme_reference}
\begin{cases}
\left(\frac{1}{\tau_n}\hat{\vec{M}}^n+D_1\hat{\vec{S}}^n+\gamma\hat{\vec{N}}_1^n\right)\hat{\vec{W}}^n_1&=\frac{1}{\tau_n}\hat{\vec{M}}^{n-1}\hat{\vec{W}}^{n-1}_1+\gamma{a}\hat{\vec{F}}^n\\
\left(\frac{1}{\tau_n}\hat{\vec{M}}^n+D_2\hat{\vec{S}}^n+\gamma\hat{\vec{N}}_2^n\right)\hat{\vec{W}}^n_2&=\frac{1}{\tau_n}\hat{\vec{M}}^{n-1}\hat{\vec{W}}^{n-1}_2+\gamma{b}\hat{\vec{F}}^n,
\end{cases}
\end{equation}
where $\vec{W}_1$ and $\vec{W}_2$ represent the nodal values of the discrete solutions corresponding to $\hat{u}_1$ and $\hat{u}_2$ respectively and the equations are nondimensional such that either $D_1$ or $D_2$ is equal to 1.
The components of the weighted mass matrix $\hat{\vec{M}}$, the weighted  stiffness matrix $\hat{\vec{S}}$ and the load vector $\hat{\vec{F}}$ on the reference frame are given by
\begin{equation}
\hat{M}^n_{\alpha\beta}:=\int_{\Oc}J^n\hat{\Phi}_\alpha\hat{\Phi}_\beta,\quad
\hat{S}^n_{\alpha\beta}:=\int_{\Oc}[J\vec{K}]^n\nabla\hat{\Phi}_\alpha\cdot\vec{K}^n\nabla\hat{\Phi}_\beta
\quad\text{and}\quad
\hat{F}^n_{\alpha}:=\int_{\Oc}J^n\hat{\Phi}_\alpha.
\end{equation}\changesvtwo{
For reaction kinetics (\ref{eqn:schnak}) the components of the matrices arising from the Picard linearisation $\hat{\vec{N}}_1$ are given by
\begin{equation}
\left(\hat{N}_1\right)_{\alpha\beta}:=\sum_{\eta=1}^{\dim(\Vc)}\sum_{\vartheta=1}^{\dim(\Vc)}\left[(W_2)_\eta(W_2)_\vartheta\right]^{n-1}\int_{\Oc}J^n\hat{\Phi}_\alpha\hat{\Phi}_\beta\hat{\Phi}_\eta\hat{\Phi}_\vartheta,
\end{equation}
with  $\hat{\vec{N}}_2$ treated similarly. 
\condense{
and
\begin{equation}
\left(\hat{N}_2\right)_{\alpha\beta}:=\int_{\Oc}J^n\hat{\Phi}_\alpha\hat{\Phi}_\beta+\sum_{\eta=1}^{\dim(\Vc)}\sum_{\vartheta=1}^{\dim(\Vc)}\left[(W_1)_\eta(W_2)_\vartheta\right]^{n-1}\int_{\Oc}J^n\hat{\Phi}_\alpha\hat{\Phi}_\beta\hat{\Phi}_\eta\hat{\Phi}_\vartheta,
\end{equation}
respectively.}}
\deletions{
By the definition of $\Vtn$ (\ref{eqn:fe_space_mov}) we obtain the following {\it time dependent} mass $(M^n)$ and stiffness $(S^n)$ matrices and load vector $\vec F^n$ if we assemble the linear systems on the evolving domain
\begin{equation}
{M}^n_{\alpha\beta}:=\int_{\Otn}\Phi^n_\alpha\Phi^n_\beta =\hat{M}^n_{\alpha\beta},\quad
{S}^n_{\alpha\beta}:=\int_{\Otn}\nabla{\Phi}^n_\alpha\cdot\nabla\Phi^n_\beta=\hat{S}^n_{\alpha\beta}\quad\text{and}\quad
{F}^n_{\alpha}:=\int_{\Otn}{\Phi}^n_\alpha=\hat{F}^n_{\alpha}.
\end{equation}
We thus obtain the following linear systems: 
\begin{equation}
\label{eqn_matrix_vector_scheme_evolving}
\begin{cases}
\left(\frac{1}{\tau_n}{\vec{M}}^n+D_1{\vec{S}}^n+\gamma{\vec{N}}_1^n\right){\vec{W}}^n_1&=\frac{1}{\tau_n}{\vec{M}}^{n-1}{\vec{W}}^{n-1}_1+\gamma{a}{\vec{F}}^n\\
\left(\frac{1}{\tau_n}{\vec{M}}^n+D_2{\vec{S}}^n+\gamma{\vec{N}}_2^n\right){\vec{W}}^n_2&=\frac{1}{\tau_n}{\vec{M}}^{n-1}{\vec{W}}^{n-1}_2+\gamma{b}{\vec{F}}^n,
\end{cases}
\end{equation}
using analogous modifications for $\vec N_1$ and $\vec{N}_2$. 
\condense{
where for reaction kinetics (\ref{eqn:schnak}) the components of the matrices arising from the Picard linearisation  $\vec{N}_1=\hat{\vec{N}}_1$ are given by
\begin{align}
\left({N}_1\right)_{\alpha\beta}&:=\sum_{\eta=1}^{\dim(\Vc)}\sum_{\vartheta=1}^{\dim(\Vc)}\left[(W_2)_\eta(W_2)_\vartheta\right]^{n-1}\int_{\Otn}{\Phi}^n_\alpha{\Phi}^n_\beta{\Phi}^n_\eta{\Phi}^n_\vartheta
\end{align}
}
}

Formulation (\ref{eqn_matrix_vector_scheme_reference}) gives rise to the following linear algebra problem: Solve for vectors $\vec{b}_i^n, i=1,\dotsc,m,$ such that
\begin{equation}
\vec{A}^n\vec{b}_i^n=\vec{c}^{n-1}_i, \text{ for }  n=1,\dotsc,N.
\end{equation}
The matrix $\vec{A}^n$ is symmetric sparse and positive definite. We therefore use the conjugate gradient  (CG) algorithm \citep{hestenesmethods}  to compute the solution to the linear systems.
\deletions{
\begin{figure}[h]
\centering
  \begin{tikzpicture}[scale=0.7] 
    \path (0,-2) node {Triangulation of $\Oc$};
    \path (7.5,-2) node {Triangulation of $\Ot$};
    \newcommand{\tikzmeshsize}{4}
    \newcommand{\arrowcushion}{-0.1}
    \coordinate (a) at (-2,-1);
    \coordinate (b) at ($(a)+(\tikzmeshsize,0)$);
    \coordinate (c) at ($(a)+(0,\tikzmeshsize)$);
    \coordinate (d) at ($(a)+(\tikzmeshsize,\tikzmeshsize)$);
    \draw[fill=b,fill opacity=0.4] (a)--(b)--(c)--cycle;
    \draw[fill=c,fill opacity=0.4] (a)--(d)--(c)--cycle;
    \draw(d)--(b); 
    \draw[stealth-stealth] ($(a)+(0,\arrowcushion)$)--($(b)+(0,\arrowcushion)$) node[pos=0.5,below] {$\hat{h}$};
       
      \draw[fill=b,fill opacity=0.4]{ 
        (4,-1) coordinate (A)
        .. controls (7.5,0) ..
        (10,-1) coordinate (B)
        .. controls (7.75,1.)  ..
        (5,4) coordinate (C)
        .. controls (5,1.5) .. 
        (4,-1)
      };
      \draw[fill=c,fill opacity=0.4]{
        (A)
        .. controls (7,.75) .. 
        (11,5) coordinate (D)
        .. controls (7.5,5) ..
        (C) 
        .. controls (5,1.5) .. 
        (A)
      };
      \draw (B) .. controls (11,1.5) .. (D);
      \draw[-stealth]{
        ($(c)!.3!(b)+(.5,.5)$) coordinate (xi) 
        --
        ($(C)!.3!(B)+(.8,1.1)$) coordinate (x)
        node[pos=0.5,above] {$\A_t$}
      };
      \draw[fill=black] (xi) circle (2pt) node[above] {$\vec{\xi}$} ;  
      \draw[fill=black] (x) circle (2pt) node[anchor=south west] {$\vec{x}=\A_t(\vec{\xi})$};  
  \end{tikzpicture}
  \caption[``Triangulation'' of the evolving domain]{A simple example of the reference and evolving domain with the associated mapping, mesh-size and triangulations.}  
  \label{Fig:mapping}   
\end{figure}
\begin{Rem}[Quadrature on the evolving domain]
\condense{
  As we do not have to compute the Jacobian of the mapping, assemblage
  of the linear systems is faster on the evolving domain. However, in
  the previous analysis we have neglected errors due to variational
  crimes, such as the fact that integrals of finite element functions
  must be evaluated by some numerical quadrature. In light of this, it
  should be noted that the finite dimensional space $\Vt$ on the
  evolving domain will not in general consist of piecewise polynomial
  functions. Furthermore, simplexes on the evolving domain will not in
  general be affine transformations of the reference simplex (see Figure
  \ref{Fig:mapping} for an example). If formulation
  (\ref{eqn:moving_fully_discrete_scheme}) is used and domain evolution
  is not spatially linear, the influence of numerical quadrature on the
  accuracy of the scheme should be considered. We leave this extension
  for future studies.
}\changesvtwo{
  In the previous analysis we have neglected errors due to variational
  crimes, such as numerical quadrature. The finite dimensional space
  $\Vt$ on the evolving domain does not generally consist of piecewise
  polynomial functions. Furthermore, simplexes on the evolving domain
  will not in general be affine transformations of the reference
  simplex (see Figure \ref{Fig:mapping} for an example). If
  formulation (\ref{eqn:moving_fully_discrete_scheme}) is used and
  domain evolution is not spatially linear, the influence of numerical
  quadrature on the accuracy of the scheme should be considered. We
  leave this extension for future studies.
}
\end{Rem}
}
\section{Numerical experiments}\label{secn.:numerics}
 We now provide numerical evidence to
back-up the estimate of Theorem \ref{thm:fully_discrete}.  We use as a
test problem, the Schnakenberg kinetics, although any other reaction
kinetics that fulfils our assumptions could have been used.
For the implementation we make use of the toolbox \ALBERTA \citep{schmidt2005design}. The graphics were generated with \PARAVIEW \citep{henderson2004paraview}. 
\condense{
\begin{Defn}[Experimental order of convergence]
We denote the $\Lp{\infty}{\left(0,T;\Lp{2}{(\Omega)^m}\right)}$ error in the numerical scheme on a series of uniform refinements of a triangulation $\big\{\hat{\T}_i\big\}_{i=0,\dotsc,N}$ by $\{ \vec{e}_i \}_{i = 0,\dotsc,N}$.
The  experimental order of convergence
(EOC) is defined to be the numerical measure of rate of convergence of the scheme as $\hat{h}_n 
\rightarrow 0$, where $\hat{h}_n$ denotes the maximum mesh-size of $\hat{\T}_n$  and is given by
\begin{equation}\label{eqn:eoc_defn}
\operatorname{EOC}_i(\vec{e}_{i,i+1},\hat{h}_{i,i+1}) = \frac{ \ln({\vec{e}_{i+1}}/{\vec{e}_i}) }{ \ln({\hat{h}^{i+1}}/{\hat{h}_i}) }.
\end{equation}
\end{Defn}
}
\changesvtwo{
\subsection{Numerical verification of the a priori convergence rate}
We examine the experimental order of convergence (EOC) of scheme (\ref{eqn:fully_discrete_scheme}).  The
EOC is a numerical measure of the rate of convergence of the scheme as $\hat{h}_n 
\rightarrow 0$. For  a series of uniform refinements of a triangulation $\big\{\hat{\T}_i\big\}_{i=0,\dotsc,N}$ we denote by $\{ \vec{e}_i \}_{i = 0,\dotsc,N}$ the error and $\hat{h}_n$  the maximum mesh-size of $\hat{\T}_n$. The EOC is given by
\begin{equation}\label{eqn:eoc_defn}
\operatorname{EOC}_i(\vec{e}_{i,i+1},\hat{h}_{i,i+1}) = \ln({\vec{e}_{i+1}}/{\vec{e}_i})/\ln({\hat{h}^{i+1}}/{\hat{h}_i}) .
\end{equation}
We  consider the EOC in approximating the solution to  (\ref{eqn:model_problem}),   with $\mathbb{P}^{1}$, $\mathbb{P}^{2}$ and $\mathbb{P}^{3}$  basis functions and uniform timestep $\tau\approx\hat{h}^2$, $\tau\approx\hat{h}^3$ and $\tau\approx\hat{h}^4$   respectively (since the scheme is first order in time). 
}
We also consider two different forms of domain evolution.
\begin{itemize}
\item{Spatially linear periodic evolution:}
\begin{equation}\label{eqn:slp_evolution}
\A_t(\vec{\xi})=\vec{\xi}\left(1+\kappa\sin\left({\pi{t}}/{T}\right)\right).
\end{equation}
\item{Spatially nonlinear periodic evolution:}
\begin{equation}\label{eqn:snp_evolution}
(\A_t(\vec{\xi}))_i={\xi}_i\left(1+\kappa\sin\left({\pi{t}}/{T}\right){\xi}_i\right)
\text{ for }i=1,\dotsc,d,
.
\end{equation} 
\end{itemize}
In both cases we take a  time interval of $[0,1]$, the initial domain as the unit square and the parameter $\kappa=1$. We take the diffusion coefficients $\vec{D}=(0.01,1)^\transpose$ and the parameter $\gamma=1$. Problem \ref{Pbm:RDS_td} equipped with nonlinear reaction kinetics  does not admit any closed form solutions. In order to provide numerical verification of the convergence rate, we insert a source term   such that the exact solution is, 
\begin{equation}\label{eqn:benchmark_exact}
\begin{split}
\hat{u}_1\left(\vec{\xi},t\right)=\sin(\pi{t})\cos(\pi x_1)\cos(\pi x_2),\quad\hat{u}_2\left(\vec{\xi},t\right)=-\sin(\pi{t})\cos(\pi x_1)\cos(\pi x_2).
\end{split}
\end{equation}

Tables \ref{tab:eoc_schnak_lin} and \ref{tab:eoc_schnak_nonlin} show the EOCs for the two benchmark examples. In both examples we observe that the error converges at the expected rate, providing numerical evidence for the estimate of Theorem \ref{thm:fully_discrete}.  
\begin{table}[htbp]
\centering
\begin{tabular}{cccccc}\toprule 
&
$\abs{\log_2\hat{h}}$
&
4
&
5
&
6
&
7
\\
\midrule
\multirow{2}{*}{$\mathbb{P}^{1}$}
&
$e$
&
2.34
e-1
&
6.20
e-2
&
1.57
e-2
&
3.96e-3
\\
&EOC
&
n.a.
&
1.91
&
1.98
&
1.99
\\
\midrule\multirow{2}{*}{$\mathbb{P}^{2}$}
&
$e$
&
3.93
e-2
&
4.96e-3
&
6.20e-4
&
8.00e-5
\\
&
EOC
&
n.a.
&
2.98
&
3.00
&
2.99
\\
\midrule\multirow{2}{*}{$\mathbb{P}^{3}$}&$e$
&
9.66e-3
&
6.10e-4
&
n.a.
&
n.a.
\\
&
EOC
&
3.89
&
3.99
&
n.a.
&
n.a.
\\
\bottomrule
\end{tabular}
\caption[Experimental order of convergence (spatially linear periodic growth)]{Error in the $\Lp{\infty}{\left(0,T;\Lp{2}{(\Oc)^m}\right)}$ norm  and EOCs for a benchmark problem  with spatially linear domain evolution (\ref{eqn:slp_evolution}). }\label{tab:eoc_schnak_lin}
\end{table}
\begin{table}
\centering
\begin{tabular}{cccccc}\toprule
  &
  $\abs{\log_2\hat{h}}$
  &
  4
  &
  5
  &
  6
  &
  7
  \\
  \midrule
  \multirow{2}{*}{$\mathbb{P}^{1}$}
  &
  $e$
  &
  1.42
  e-1
  &
  3.63
  e-2
  &
  9.12
  e-3
  &
  2.28
  e-3
  \\
  &
  EOC
  &
  n.a.
  &
  1.97
  &
  1.99
  &
  1.99
  \\
  \midrule
  \multirow{2}{*}{$\mathbb{P}^{2}$}
  &
  $e$
  &
  1.07
  e-2
  &
  1.32
  e-3
  &
  1.60
  e-4
  &
  2.00
  e-5
  \\
  &
  EOC
  &
  n.a.
  &
  3.02
  &
  3.02
  &
  3.01
  \\
  \midrule
  \multirow{2}{*}{$\mathbb{P}^{3}$}
  &
  $e$
  &
  1.89
  e-3
  &
  1.20
  e-4
  &
  n.a.
  &
  n.a.
  \\
  &
  EOC
  &
  3.97
  &
  4.00
  &
  n.a.
  &
  n.a.
  \\
  \bottomrule
\end{tabular}
\caption[Experimental order of convergence (spatially nonlinear periodic growth)]{Error in the $\Lp{\infty}{\left(0,T;\Lp{2}{(\Oc)^m}\right)}$ norm  and EOCs for a benchmark problem with  nonlinear domain evolution (\ref{eqn:snp_evolution}).  }\label{tab:eoc_schnak_nonlin}
\end{table}

\additions{
\begin{Rem}[Existence of solutions to Problem \ref{Pbm:RDS_td} with spatially linear isotropic evolution]
In \cite{venkataraman2010global}, we showed that Problem
\ref{Pbm:RDS_td} equipped with the Schnakenberg reaction kinetics posed on
a $C^2$ domain $\Ot$, is well posed under any bounded spatially linear
isotropic evolution of the domain.  If we assume this result holds on
polygonal domains, we have sufficient regularity on the continuous
problem to apply Theorem \ref{thm:fully_discrete} and thus conclude
scheme (\ref{eqn:fully_discrete_scheme}) with $\mathbb{P}^1$ finite
elements converges with optimal order.
\end{Rem}
Thus, to illustrate a concrete application for which our theory holds, we present results for the Schnakenberg kinetics with domain growth function of the form (\ref{eqn:slp_evolution}), initial conditions are taken as small perturbations around the spatially homogeneous steady state and numerical and reaction kinetic parameter values as given in Table \ref{tab:parameter_values}.
\begin{table}[h]
\centering
\begin{tabular}{  cccccccccc }
\toprule
$D_1$&$D_2$&$\gamma$&$a$&$b$&$\kappa$&$T$&$\tau$&DOFs\\
.01&1.0&0.1&0.1&0.9&4&2000&$10^{-2}$&8321\\
\bottomrule
\end{tabular}\\
\caption[]{Parameter values for the numerical experiment with the Schnakenberg kinetics.}\label{tab:parameter_values}
\end{table}

We take the unit square as the initial domain, with the domain growing from a square of length 1 to a square of length 5 at $t=1000$ before contracting to a square of length 1 at final time. Figure \ref{fig:lagrangian-patterns} shows snapshots of the discrete activator ($W_1$) profiles. The substrate profiles ($W_2$) have been omitted as they are $180^\circ$ out of phase with those of the activator. An initial half spot pattern forms which reorients as the domain grows into a single spot positioned in the centre of the domain. As the domain contracts this single spot disappears (via spot annihilation)  with the final domain exhibiting no spatial patterning. 
\begin{figure}
\centering
\includegraphics[trim = 0mm 0mm 0mm 0mm, clip=true, width=0.75\linewidth]{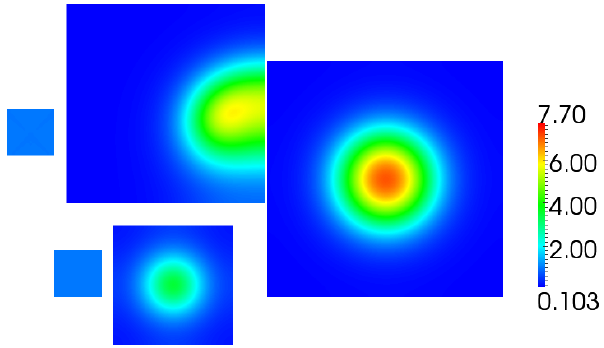}
\caption{Snapshots of the discrete activator ($u_1$) profile  for the Schnakenberg reaction kinetics on domains with spatially linear evolution at times 0, 590, 1000, 1750 and 2000 reading clockwise from top left.  For  parameter values see  Table \ref{tab:parameter_values}. We observe the formation of a half spot which reorients to a single spot positioned in the centre of the domain. As the domain contracts the spots are annihilated with the domain at end time exhibiting no patterns. }\label{fig:lagrangian-patterns}  
\end{figure}
}
\deletions{
We now present numerical results illustrating the influence of domain evolution on pattern formation by RDSs.
We first present results for the Schnakenberg kinetics with domain growth functions (\ref{eqn:slp_evolution}) and (\ref{eqn:snp_evolution}) with the same initial conditions and numerical and reaction kinetic parameter values as given in Table \ref{tab:parameter_values}.
\begin{table}[h]
\centering
\begin{tabular}{  cccccccccc }
\toprule
$D_1$&$D_2$&$\gamma$&$a$&$b$$\kappa$&$T$&$\tau$&DOFs\\
.01&1.0&0.1&0.1&0.94&2000&$10^{-2}$&8321\\
\bottomrule
\end{tabular}\\
\caption[]{Parameter values for numerical experiments with the Schnakenberg kinetics.}\label{tab:parameter_values}
\end{table}

We take the unit square as the initial domain and the evolution of the boundary curve of the domain is identical  in both cases, with the domain growing from a square of length 1 to a square of length 5 at $t=1000$ before contracting to a square of length 1 at final time. Figure \ref{fig:lagrangian-patterns} shows snapshots of the discrete activator ($W_1$) profiles obtained using the different schemes. The substrate profiles ($W_2$) have been omitted as they are $180^\circ$ out of phase with those of the activator. In the spatially linear case (Figure \ref{fig:lagrangian-patterns} left), an initial half spot pattern forms which continuously transitions as the domain grows into a single spot positioned in the centre of the domain. As the domain contracts this single spot disappears (via spot annihilation)  with the final domain exhibiting no spatial patterning. In the spatially nonlinear growth case (Figure \ref{fig:lagrangian-patterns}, right) the pattern transition is completely different with a half spot forming which splits as the domain grows to form two half spots which move around the domain before being annihilated as the domain contracts. The difference in patterning observed appears to be due to the differences in the growth  function, with the results illustrating the robustness of the  numerical method in dealing with non-uniform domain evolution that is likely to be encountered in the  biological problems we have in mind.

 Along with 2 component RDSs, 3 component RDSs where a second
 inhibitor quenches established maximums have been widely studied
 \citep{meinhardt2004out}. Many interesting phenomena, which do not
 occur in 2 component systems, such as out of phase oscillations and
 space-time patterning which does not reach a steady state (even on
 fixed domains) are observed. The applications of such 3 component
 systems are of much importance, for example in the modelling of cell
 polarisation during chemotaxis \citep{meinhardt1999orientation}.
 \condense{
   To illustrate the versatility of our method in dealing with
   multiple component systems, we now present results for a three
   component model. We first consider a fixed domain, to illustrate
   both the qualitatively different solution behaviour to the 2
   component case and that our method is readily applicable to
   multicomponent RDSs posed on fixed domains.
 }
 \begin{Defn}[Global and local inhibition kinetics \citep{meinhardt2004out}]
   The following model is comprised of a single self-enhancing activator $u_1$ antagonised by global (fast diffusing) and local (slow diffusing)  inhibitors; $u_2$ and $u_3$ respectively:
   \begin{equation}\label{eqn:3component}
     \begin{split}
         f_1\left(u_{1}, u_{2}, u_{3}\right)
         &
         =
         \gamma\left(
         \frac{s\left(u_{1}^2+b_{1}\right)}{u_2\left(1+s_3u_3\right)}-r_1u_1
         \right),
         \\
         f_2\left(u_{1}, u_{2}, u_{3}\right)
         &
         =
         \gamma\left(s u_1^2-r_2u_2\right),
         \\
         f_3\left(u_{1}, u_{2}, u_{3}\right)
         &
         =
         \gamma\left(r_3u_1-r_3u_3\right),
     \end{split}
   \end{equation}
   where $0 \leq s, b_1, r_1,  s_3, r_2,  r_3, \gamma < \infty $.
\end{Defn}
Figure \ref{fig:3component_evolving} shows snapshots of the numerical
solution on a domain with growth of the form (\ref{eqn:slp_evolution})
and parameter values as given in Table
\ref{tab:3component_parameter_values}.  Broadly speaking similar
behaviour to the 2 component case is observed with the patterning mode
(number of spots) generally increasing as the domain grows and
decreasing as the domain contracts. The initial pattern is a 2 peak
pattern that exhibits travelling wave behaviour similar to the fixed
case. After the domain has grown sufficiently large new peaks appear
either via splitting of existing peaks or peak insertion. When
multiple peak patterns are present the behaviour is much more
complicated than the two component case, the peaks travel around the
domain and can collide leading to peak-merging or when peaks are in
close proximity peak-annihilation. This behaviour of peak-merging and
annihilation was observed previously by \citet{venkataraman2010global}
for 2 component systems, the novelty in the 3 component case is that
this behaviour is observed even as the domain grows. This is of
interest as with this specific 3 component system domain growth
(contraction) does not lead to a monotonic increase (decrease) in the
number of peaks.
\begin{table}
\centering
\begin{tabular}{  cccccccccccccc }
\toprule
$D_1$&$D_2$&$D_3$&$\gamma$&$s$&$s_3$&$r_1$&$r_2$&$r_3$&$b_1$&$\kappa$&$T$&$\tau$&DOFs\\
1.0&100&0.01&$10^{4}$&$0.005$&0.8&0.005&0.008&0.001&0.01&4&100&$10^{-3}$&8321\\
\bottomrule
\end{tabular}\\
\caption[]{Parameter values for numerical experiments with the 3 component kinetics.}\label{tab:3component_parameter_values}
\end{table}
\begin{figure}
\centering
\includegraphics[trim = 0mm 220mm 0mm 0mm, clip=true, width=\linewidth]{./figures/schnakenberg.pdf}
\caption{Snapshots of the discrete activator ($u_1$) profile  for the Schnakenberg reaction kinetics on domains with spatially linear (left) and nonlinear (right) evolution at times 0, 590, 1000, 1750 and 2000 reading clockwise from top left.  For  parameter values see  Table \ref{tab:parameter_values}. Although the evolution of the domain boundary is identical in both cases, the pattern transitions observed are markedly different. Under spatially linear growth we observe the formation of a half spot which reorients to a single spot positioned in the centre of the domain. Under nonlinear growth, we observe the formation of the same initial half spot that now splits into two half spots. In both cases,  as the domain contracts the spots are annihilated with the domain at end time exhibiting no patterns. }\label{fig:lagrangian-patterns}  
\end{figure}
\begin{figure}
\centering
\includegraphics[trim = 35mm 210mm 35mm 19mm, clip=true, width=0.49\linewidth]{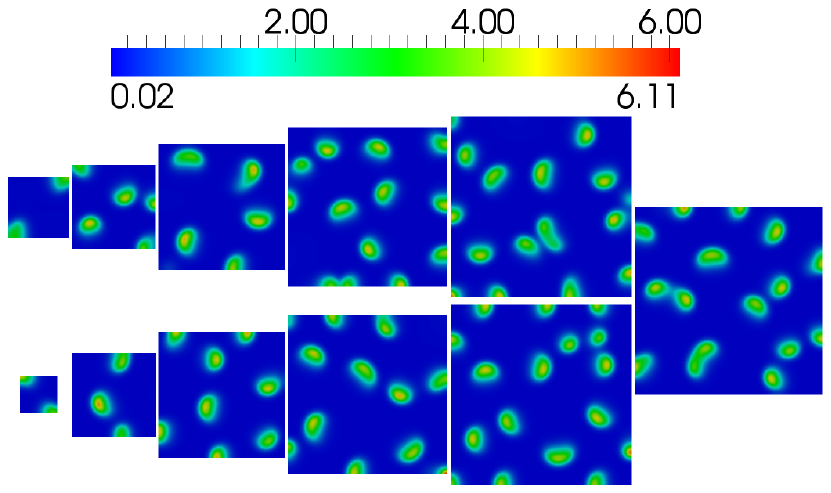}
\includegraphics[trim = 35mm 210mm 35mm 19mm, clip=true, width=0.49\linewidth]{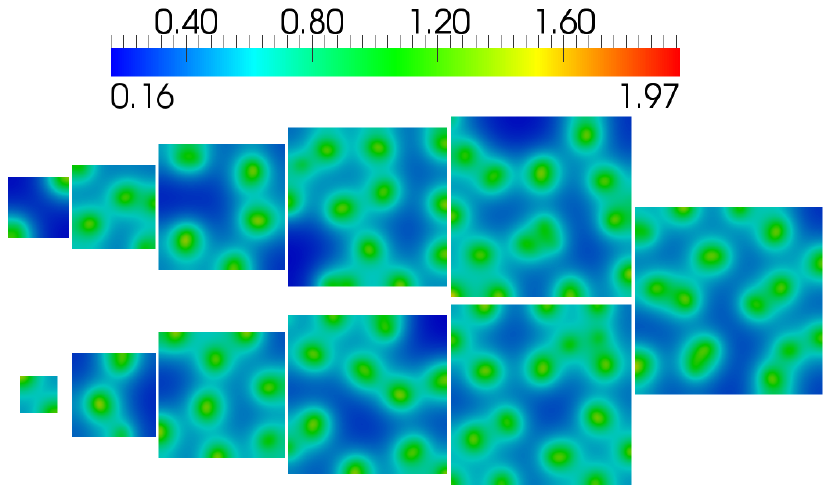}\\
\includegraphics[trim = 35mm 210mm 35mm 19mm, clip=true, width=0.49\linewidth]{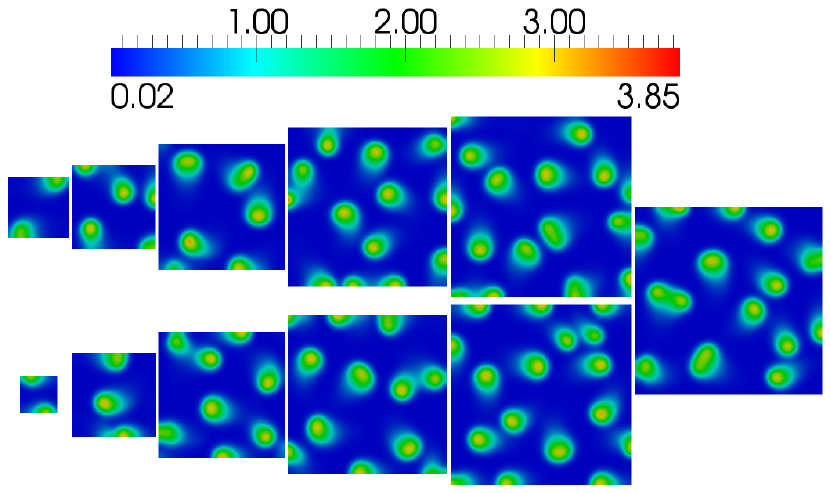}
\caption{Discrete activator $u_1$ (top left), global inhibitor $u_2$ (top right) and local inhibitor $u_3$ (bottom) concentrations of the 3 component system (\ref{eqn:3component}) on a domain with spatially linear evolution at times 5, 10,  20, 30, 40, 50, 60, 70, 80, 90 and 100 reading clockwise from top left.  The colour (grey-scale) legend indicates the numerical range of the solution during the experiment.  Generally more peaks are observed on larger domains, we observe  travelling wave like solutions and peak collisions which result in spot merging as well as peak annihilation.}\label{fig:3component_evolving}  
\end{figure}
}
\bibliographystyle{abbrvnat}
              
\end{document}